\documentclass[11pt]{article}
\usepackage{amsmath,amsfonts,amssymb,mathrsfs,tikz-cd,dsfont,amsthm,tikz,cjhebrew,authblk,titling, hyperref}
\usepackage[margin=1in]{geometry}
\usepackage[backend=bibtex,eprint=true, url=false,doi=false,isbn=false,style=alphabetic]{biblatex}
\bibliography{bibliography}
\numberwithin{equation}{subsection}
\hypersetup{colorlinks=true,linkcolor=blue,citecolor=black,linktoc=all}

\newtheorem{theorem}{Theorem}[subsection]
\newtheorem{lemma}[theorem]{Lemma}
\newtheorem{proposition}[theorem]{Proposition}
\newtheorem{corollary}[theorem]{Corollary}

\theoremstyle{definition}
\newtheorem{definition}[theorem]{Definition}
\newtheorem{remark}[theorem]{Remark}
\newtheorem*{proof*}{proof}
\newtheorem{example}[theorem]{Example}
\newtheorem{notation}[theorem]{Notation}
\newtheorem{construction}[theorem]{Construction}
\newtheorem{warning}[theorem]{Warning}
\newtheorem{conjecture}{Conjecture}
\newtheorem{hypothesis}[theorem]{Hypothesis}

\newcommand{\End}{End}
\newcommand{\LM}{\mathrm{LM}}
\newcommand{\id}{\mathrm{id}}
\newcommand{\lax}{\mathrm{lax}}
\newcommand{\op}{\mathrm{op}}
\newcommand{\rB}{{\mathrm{B}}}
\newcommand{\psh}{\mathrm{PSh}}
\newcommand{\opd}{\mathrm{Opd}}
\newcommand{\fun}{\mathrm{Fun}}

\newcommand{\cat}{\mathrm{Cat}}

\newcommand{\Vect}{Vect}

\newcommand{\Lmod}{\mathrm{LMod}}
\newcommand{\Rmod}{\mathrm{RMod}}
\newcommand{\bimod}{\mathrm{Bimod}}
\newcommand{\Mod}{\mathrm{Mod}}
\newcommand{\alg}{\mathrm{Alg}}
\newcommand{\EM}{\mathrm{EM}}

\newcommand{\colim}{\mathrm{colim}}
\let\lim\relax
\newcommand{\lim}{\mathrm{lim}}
\let\pr\relax
\newcommand{\pr}{\mathrm{Pr}}

\DeclareFontFamily{U}{min}{}
\DeclareFontShape{U}{min}{m}{n}{<-> udmj30}{}
\newcommand\yo{\!\text{\usefont{U}{min}{m}{n}\symbol{'207}}\!}
\newcommand{\FZ}{\text{\usefont{U}{euf}{m}{n}Z}}

\DeclareFontFamily{U}{rcjhbltx}{}	
\DeclareFontShape{U}{rcjhbltx}{m}{n}{<->rcjhbltx}{}
\DeclareSymbolFont{hebrewletters}{U}{rcjhbltx}{m}{n}

\DeclareMathSymbol{\tav}{\mathord}{hebrewletters}{116}

\title{Large condensation in enriched $\infty$-categories}
\author{Devon Stockall\thanks{stockall@imada.sdu.dk}}
\date{\vspace{-5ex}}
\begin{document}
\maketitle

\begin{abstract}
We demonstrate that many algebraic features of higher condensation theory are really features of higher algebra.  We work in the context of monads on the unit of a monoidal $(\infty,2)$- category.  Assuming that 2-categorical unstraightening is monoidal, we show that the functor which takes such a monad to its Eilenberg-Moore object is lax monoidal.  We use this to introduce an analog of condensation for codimension-i symmetries, without assuming dualizability.  We give sufficient conditions for Morita equivalence of endomorphisms of condensates. We relate higher centralizers to iterated categories of modules.  Along the way, we prove results pertaining to functoriality of Day convolution for enriched $\infty$-categories, which may be of independent interest.

\end{abstract}

\tableofcontents

\section{Introduction}

The categorical structure of three-dimensional anyon condensation (or 1-form symmetry gauging) is justified in \cite{Kong2014}.  One begins with the data of a modular tensor category (MTC) $\mathcal{C}$, which describes line operators in a three-dimensional TQFT, and an $E_2$-algebra object $A$ in $\mathcal{C}$ corresponding to condensable anyons. The category $\Mod_A^{\mathrm{loc}}(\mathcal{C})$ of \emph{local $A$-modules} in $\mathcal{C}$ is found to describe line operators in the theory in which $A$ is condensed.  This theory relies heavily on semisimplicity, dualizability, and seperability, which is encapsulated in \emph{full dualizability} \cite{DSPS}.  By assuming full-dualizability and making use of the Baez-Dolan cobordism hypothesis, the fusion-categorical theory of condensation has been extended to arbitrary dimension \cite{BaezDolan,LurieTFT,KongBulkBoundary,KongBulkBoundaryII,GaiottoJohnsonFreyd,DouglasReutter2018,JohnsonFreyd,QuantumLiquids,Kong2024}.  

Another similar construction in two-dimensional conformal field theory (CFT) describes representations of a vertex operator algebra (VOA) in terms of representations of a vertex subalgebra \cite{HKL,CKM}.  Given a VOA $V$, an abelian braided category $\mathcal{C}$ of $V$-representations, and VOA extension $A\supset V$, one considers $A$ as an $E_2$-algebra in $\mathcal{C}$.  The category of $A$-representations, as a VOA, is braided equivalent to $\Mod_A^{\mathrm{loc}}(\mathcal{C})$.  To address logarithmic CFTs, the theory of VOA extension was developed to function without requiring the full-dualizability that features in the fusion setting. 

Our goal is to show that many algebraic features of condensation theory are really features of higher algebra.  We prove mathematical results, and then use them in Section \ref{condensationsec} to provide common generalization of these two examples, allowing condensation of symmetries in arbitrary (co)dimension, without requiring dualizability.  The introduction of higher-algebraic tools was likely unavoidable in the construction of this generalization: algebraic theories of $n$-categories quickly become intractable in higher dimensions.  We present this theory in the language of enriched $\infty$-categories, for which many fundamental results have been proven only recently \cite{GepnerHaugseng,Hinich2021,Hinich2023,HeineEquivalence,HeineBienriched, HeineWeighted, Ben-Moshe2024,ReutterZetto}.  This allows additional flexibility.  Notably, $\infty$-categories `remember' homotopy-type by nature, allowing for homotopical condensation of continuous symmetries, previously inaccessible to the categorical theory.

\subsection{Results}
For $2$-category $\mathfrak{D}$, a \emph{monad} on $d\in \mathfrak{D}$ is an algebra object in the $1$-category of endomorphisms $\Omega_d \mathfrak{D}:= \mathfrak{D}(d,d)$.  To a monad, we can assign an Eilenberg-Moore object, which is a `universal module' over this monad (assuming such a universal object exists).  The physical role of such an object is discussed in Example \ref{EM1}.  This will be the primary technical tool of this work.  The Eilenberg-Moore object $\EM(T)$ of a monad $T$ on $d\in \mathfrak{D}$ is equipped with an adjunction 
\begin{equation}
	\begin{tikzcd}
		d\ar[r,shift left=2,"\EM_L(T)"{name=A}]& \EM(T)\ar[l, shift left=2, "\EM_R(T)"{name=B}]\ar[from=A,to=B, phantom,"\dashv"rotate=-90]
	\end{tikzcd}
\end{equation}
The assignment of a monad to its Eilenberg-Moore object can be made functorial in two ways: the \emph{right Eilenberg-Moore functor} $T\mapsto \EM_R(T)$, or, if $\mathcal{D}$ admits the appropriate colimits, the \emph{left Eilenberg-Moore functor} $T\mapsto \EM_L(T)$.  Functors of this form are constructed in \cite{HaugsengMonad,HeineMonadicity,StefanichPhD}.  If $\mathfrak{D}$ is monoidal and $d=\mathds{1}_ \mathfrak{D}$ is the monoidal unit, then we can ask that a monad has additional algebraic structure.  The goal of this work is to prove that, if $T\in \alg(\Omega_{\mathds{1}_ \mathfrak{D}}\mathfrak{D})$ has additional algebraic structure, then its Eilenberg-Moore object does also.  
\begin{warning}
	The proofs of Theorem \ref{mainthm:EM} and Corollary \ref{MainCor:Alg} provided in this work make use of the hypothesis that 2-categorical unstraightening respects monoidal structure (Hypothesis \ref{hyp:MonoidalCoCart}).  The $1$-categorical analogue of this is proven in \cite[]{Hinich2015,RamziUnstraightening}, and the 2-categorical result is certainly also true.  Similar results can be obtained without this hypothesis using \cite[]{HeineMonadicity}. 
\end{warning}

Denote by $\mathcal{D}^R_{/d}\subset \mathcal{D}_{/d}$ the full subcategory on the right adjoint morphisms to $d$, and similarly by $\mathcal{D}^L_{d/}\subset \mathcal{D}_{d/}$ the full subcategory on the left adjoints from $d$. 
\begin{theorem}[\ref{EMR},\ref{EML}]\label{mainthm:EM}
	Suppose that $\mathcal{D}$ is an $\mathcal{O}$-monoidal $(\infty,2)$-category that admits Eilenberg-Moore objects at the monoidal unit $\mathds{1}_ \mathfrak{D}\in \mathfrak{D}$.  
	\begin{itemize}
		\item The right Eilenberg-Moore functor $\EM_R$ is fully faithful and lax $\mathcal{O}$-monoidal with left adjoint
	\begin{equation}
		\begin{tikzcd}
			\alg(\Omega\mathfrak{D})^{op}\ar[r,shift left=2,"\EM_R"{name=B}]&		\mathfrak{D}_{/\mathds{1}_ \mathfrak{D}}^R \ar[l,shift left=2,"\End_L"{name=A}] \ar[phantom,from=A,to=B,"\dashv"rotate=90] , 
		\end{tikzcd}
	\end{equation}
\item If $\mathfrak{D}$ is locally compatible with simplicial colimits, then the left Eilenberg-Moore functor is lax $\mathcal{O}$-monoidal with right adjoint
\begin{equation}
		\begin{tikzcd}
			\alg(\Omega\mathfrak{D})\ar[r,shift left=2,"\EM_L"{name=B}]&		\mathfrak{D}_{\mathds{1}_ \mathfrak{D}/}^L \ar[l,shift left=2,""{name=A}] \ar[phantom,from=A,to=B,"\dashv"rotate=-90] , 
		\end{tikzcd}
	\end{equation}
	\end{itemize}
\end{theorem}

We can use this to obtain an enriched analogue of \cite{LurieHA} section 4.8.5, and a vast generalization of \cite{Kong2024} Lemma 5.2.1, allowing us to assign, for each $E_i$-algebra object, an $E_{i-1}$-monoidal module category.  This is a single step in the iterative condensation procedure of \cite{Kong2024}. 

Let $\mathcal{V}$ be a presentably $E_{k+1}$-monoidal category, and take $\cat^\mathcal{V}(\mathcal{K})$ to be the $E_k$-monoidal category of $\mathcal{V}$-enriched categories which admit $\mathcal{K}$-indexed colimits.  Let $\mathcal{C}\in \alg_{E_k} \cat^ \mathcal{V}(\mathcal{K})$ be an $E_k$-monoidal enriched category.  We apply the previous theorem to $\mathfrak{D}=\Rmod_ \mathcal{C}(\cat^ \mathcal{V}(\mathcal{K}))$.  Then $\Omega \mathfrak{D}\simeq \mathcal{C}$.

\begin{corollary}[\ref{AlgebraEMR}, \ref{AlgebraEML}, \ref{presentableAlgEML}]\label{MainCor:Alg}
	Let $\mathcal{C}\in \alg_{E_k}\cat^\mathcal{V}(\mathcal{K})$ be an $E_k$-monoidal $\mathcal{V}$-enriched category
	\begin{itemize}
		\item The contravariant functor taking an algebra to its (pointed) category of modules 
	\begin{equation}
		\Lmod^*_{-}(\mathcal{C}):\alg(\mathcal{C})^{\mathrm{op}}\to \Rmod_\mathcal{C}(\cat^\mathcal{V}(\mathcal{K}))_{/\mathcal{C}}. 
	\end{equation}
is lax $E_{k}$-monoidal and fully faithful with left adjoint. 
		\item If $\Delta^\op\in \mathcal{K}$, then the covariant functor taking an algebra to its category of modules 
	\begin{equation}
		\Lmod_{!-}(\mathcal{C}):\alg(\mathcal{C})\to \Rmod_\mathcal{C}(\cat^\mathcal{V})_{\mathcal{C}/}.
	\end{equation}
	is lax $E_k$-monoidal with right adjoint.  Moreover, if $\mathcal{C}$ is a presentably monoidal $\mathcal{V}$-category, then this functor also fully faithful.  
\end{itemize}
\end{corollary}
We recover an analogue of the fact that, in fusion categories, condensation descendants are Morita equivalent.  This generalizes results from the finite tensor category theory of \cite{EGNO}.  
\begin{proposition}\ref{prop:MonadMorita}
	If $M\in \Rmod_A \mathcal{W}$ is right dualizable and $M$ generates $\Rmod_A \mathcal{W}$, then  
	\begin{equation}
\Rmod_A\mathcal{W}(M,-)^\mathrm{enh}:\Rmod_A \mathcal{W}\to \Rmod_{\Omega_M \Rmod_A \mathcal{W}} \mathcal{W} 
	\end{equation}
	 is an equivalence.  
\end{proposition}

Finally, we prove that higher centralizers may be computed instead as loopings of lower centralizers of categories of modules. 
\begin{proposition}[\ref{iterativecentralizer}]
	Suppose that $\mathcal{W}\in Alg_{E_{k+1}}(\pr)$, and $f:C\to D$ is a morphism in $\alg_{E_k} \mathcal{W}$.  Then one has an equivalence in $\alg_{E_k} \mathcal{W}$
\begin{equation}
	\FZ_{E_k}(f)\simeq \Omega \FZ_{E_{k-1}}(\Rmod_f \mathcal{W}). 
\end{equation}
\end{proposition}
We also use a conjecture to provide alternative means to compute centralizers as endomorphism objects in categories of operadic modules.  

\subsection{Future work}
Some constructions of fusion-categorical condensation theory are notably absent from the present work.  We would like to recreate and extend these within the current framework.  Most notably, we would like to explore operadic modules \cite[Section 3.3.3]{LurieHA} in this context.  We believe that an analogue of \cite[Theorem 5.2.8]{Kong2024} holds.  Namely, in Section \ref{centralizersec}, we conjecture that if $A$ is an $E_{k+1}$-algebra in a sufficiently nice category, then there is an equivalence 
\begin{equation}\label{eq:OperadicModule}
	\Mod_A^{E_{k+1}}\simeq \Omega \Mod^{E_k}_{\Rmod_A}.  
\end{equation}

This equivalence in dimension 3 for $k=2$ gives a \emph{non-invertible} analogue of $G$-crossed braiding for \textbf{all} VOA extensions: if $A$ is a `finite group-like' extension of $V$, then the category $\Mod_A(\mathcal{C})$ of \textit{all} $A$-modules is  $G$-crossed braided \cite{KirillovMTCs,McRaeGCrossed}.  By applying \cite{JPR2023}, this determines a monoidal $2$-category whose looping is the category of local modules.  But local modules are exactly $E_2$-modules.  By applying Equation \ref{eq:OperadicModule}, we conclude that, for any vertex algebra extension, there is a monoidal 2-category whose looping is the category of local modules, namely $\Mod^{E_1}_{\Rmod_A}\simeq \bimod_{\Rmod_A|\Rmod_A}$, and that VOA extension fits into the `iterative condensation' picture.  We propose that this is the correct notion of \emph{non-invertible crossed-braided extension} of $\Mod^{\mathrm{loc}}_A$.

\subsection{Structure of the paper}
Section \ref{Backgroundsec} contains conventions and background on $\infty$-operads and enriched $\infty$-categories.  It contains few proofs, but summarizes key results, allowing us to work essentially `synthetically', without reference to any specific model of enriched $\infty$-categories. We refer to \cite{ReutterZetto} for proofs and further details. 

Section \ref{condensationsec} applies the Eilenberg-Moore functor to develop an algebraic theory analogous to the theories of higher condensation presented in \cite{Kong2014,BhardwajTachikawa2017,KongBulkBoundary,KongBulkBoundaryII,GaiottoJohnsonFreyd,JohnsonFreyd,QuantumLiquids,Kong2024}.  In subsection \ref{centralizersec}, we prove a relationship between higher centralizers and iterated module categories.  In subsection \ref{MoritaSec}, we prove that the Morita dual of a dualizable module is Morita equivalent to the starting algebra. 

Section \ref{Monadssec} contains most of the technical proofs of this work.  It serves to prove results for functoriality of enriched Day convolution, and uses them to construct both the right and left monoidal Eilenberg-Moore functors.  We then use these to construct restriction- and extension- of scalars versions of monoidal functors that take an algebra object to its category of (enriched) modules. 

\subsection{Conventions}
For us, a \emph{category} refers to an $(\infty,1)$-category: a category with $i$-morphisms for all $i\in \mathbb{Z}_{\geq 0}$, and which are all invertible for $i>1$.  Similarly, an $n$-category refers to an $(\infty,n)$-category.  If $\mathcal{C}$ is a category, we shall write $x\in \mathcal{C}$ to indicate that $x$ is an object of $\mathcal{C}$.  For $x,y\in \mathcal{C}$, we shall denote the space of morphisms from $x$ to $y$ by $\mathcal{C}(x,y)$. 

We fix inaccessible cardinals $\tav<\widehat{\tav}$.  We denote the universe of $\tav$-small categories by $\cat$, and call its objects \emph{small categories}.  Similarly, we denote the category of $\widehat{\tav}$-small categories by $\widehat{\cat}$, and call its objects \emph{large categories}.  Let $\mathcal{S}\subset \cat$ be the full subcategory on the \emph{spaces/ groupoids/ homotopy types/ anima}: these are categories for which all morphisms are invertible.  

Let $\pr$ be the category of \emph{presentable categories}, as defined in \cite[Chapter 5]{LurieHTT}.  The objects of $\pr$ are categories which admit all colimits indexed by small categories, and are generated by a small category in an appropriate sense.  The morphisms in $\pr$ are functors which preserve small colimits.  $\pr$ is closed symmetric monoidal; its product is called the \emph{Lurie tensor product} \cite[Proposition 4.8.1.15]{LurieHA}.

\subsection{Acknowledgements}
We would like to thank Markus Zetto for advice on some of the technical results of this article, and for reading and commenting on drafts while it was in preparation.  The author would like to thank Liang Kong and the anonymous referees for detailed feedback on the first version of this work.  This work was supported by VILLUM FONDEN, VILLUM Young Investigator grant 421250.

\section{Categorical background}\label{Backgroundsec}
\subsection{Operads}

\begin{definition}\hfill 
\begin{itemize}
	\item 	We define an \emph{operad} as in \cite[Definition 2.1.1.10]{LurieHA}.  Let $\opd$ denote the category of operads. 
	\item $\opd$ is presentably closed monoidal with the \emph{Boardman-Vogt} tensor product (\cite[]{BarkanSteinebrunner}, \cite[Section 2.2.5]{LurieHA}).  That is, if $\mathcal{P},\mathcal{Q},\mathcal{R}$ are operads, then the Boardman-Vogt tensor product $\mathcal{P}\otimes \mathcal{Q}$ is a corepresenting object 
	\begin{equation}
		\opd(\mathcal{P},\opd(\mathcal{Q}, \mathcal{R}))\simeq \opd(\mathcal{P}\otimes \mathcal{Q}, \mathcal{R}). 
	\end{equation}
\item An \emph{$\mathcal{O}$-monoidal category} is an object $\mathcal{C}\in \opd_{/\mathcal{O}}$ which is also a coCartesian fibration. 
\end{itemize}
\end{definition}

\begin{example}
	For $k\in \mathbb{Z}_{\geq 0}$, the \emph{little $k$-cubes} operad $E_k$ \cite[Definition 5.1.0.2]{LurieHA} encodes $k$ coherently commuting multiplication operations using embeddings of $k$-cubes into a larger $k$-cube.  An $E_1$-monoidal category is simply a monoidal category.  An $E_2$-monoidal category is a braided monoidal category. 
\end{example}

\begin{example}
	There are operads $\LM, \mathrm{RM},\mathrm{BM}$ \cite[Definition 4.2.1.1, 4.3.1.1]{LurieHA} which encode the structure of left modules, right modules, and bimodules, respectively.  An $\LM$-monoidal category is a pair of a monoidal category and a left module over it. 
\end{example}

\begin{notation}
	Suppose that $\mathcal{P},\mathcal{Q}\in \opd_{/\mathcal{O}}$.  We shall also denote the internal hom by $\alg_{\mathcal{P}/\mathcal{O}}(\mathcal{Q}):= \opd_{/\mathcal{O}}(\mathcal{P},\mathcal{Q})$.  When $\mathcal{P},\mathcal{Q}$ are $\mathcal{O}$-monoidal categories, the objects $F\in \alg_{\mathcal{P}/\mathcal{O}}(\mathcal{Q})$ will be called \emph{lax $\mathcal{O}$-monoidal functors} or \emph{$\mathcal{P}$-algebra objects}. 
\end{notation}

\begin{proposition}[\cite{Hinich2023} Proposition 8.1.3]\label{prop:AlgMonoidal}
	Let $\mathcal{M}$ be an $\mathcal{O}\otimes \LM$-monoidal cateogry which is compatible with geometric realizations.  Then the forgetful functor 
	\begin{equation}
		\alg_{\LM/ \mathcal{O}\otimes \LM}(\mathcal{M})\to \alg_{E_1/\mathcal{O}\otimes \LM}(\mathcal{M})
	\end{equation}
is an $\mathcal{O}$-monoidal coCartesian fibration. 

\end{proposition}

\subsection{Enriched categories}

\begin{definition}
Suppose that $\mathcal{V}$ is a monoidal category and $\mathcal{M}$ is a left $\mathcal{V}$-module category, with the tensoring denoted by $\otimes:\mathcal{V}\times \mathcal{M}\to \mathcal{M}$.  
\begin{itemize}
	\item A \emph{morphism object} for $m,m'\in \mathcal{M}$ is a representing object $\mathcal{M}(m,m')\in \mathcal{V}$ 
\begin{equation}
	\mathcal{V}(-\otimes m,m')\simeq \mathcal{V}(-,\mathcal{M}(m,m')):\mathcal{V}^\op\to \mathcal{S}. 
\end{equation}
\item $\mathcal{M}$ is called \emph{closed} if a morphism object exists for each $m,m'\in \mathcal{M}$. 
\end{itemize}

\end{definition}
\noindent If $\mathcal{M}$ is closed, then it is immediate from the definition that one obtains a functor 
\begin{equation}
	\mathcal{M}(-,-):\mathcal{M}^\op\times \mathcal{M}\to \mathcal{V}. 
\end{equation}
That is, every closed $\mathcal{V}$-module category is canonically \emph{$\mathcal{V}$-enriched}.  However, it is not the case that every $\mathcal{V}$-enriched category is a module category: this requires the existence of colimits called \emph{tensors}.  

\begin{notation}
	$\mathcal{V}\in \alg(\pr)$ is called \emph{presentably monoidal}. Denote $\pr^ \mathcal{V}:= \Lmod_ \mathcal{V}(\pr)$. \end{notation}
If $\mathcal{M}\in \pr^\mathcal{V}$, then it is automatically a closed module category by adjoint functor theorem, so it is automatically $\mathcal{V}$-enriched.  Conversely, all tensors exist.  This is to say, if we assume the existence of enough (co)limits, then the notions of enrichment and tensoring coincide.  This was proven in \cite{HeineEquivalence} that this is an equivalence of categories. 

More generally, one can think of a $\mathcal{V}$-enriched category as a sort of `weak left module' $\mathcal{M}$ for $\mathcal{V}$, equipped with a closed structure, giving a hom object for each pair of objects in $\mathcal{M}$. This perspective is formalized in \cite{LurieHA} 4.2.1.28.  However, this definition is technically difficult to work with.  Alternative definitions of $\mathcal{V}$-enriched categories appear in \cite{GepnerHaugseng,Hinich2021,ReutterZetto}, and equivalence between the various models is proven in \cite{HeineEquivalence,Macpherson2019}.  Enriched versions of the Yoneda embedding are given in \cite{Hinich2021,HeineEquivalence,ReutterZetto} and are proven to be equivalent in \cite{HeinePresheaves}.  Further properties are explored in \cite{Hinich2023,Ben-Moshe2024,HeineBienriched,HeineWeighted}.

We will not give the definition or provide proofs here, since we will not need to work with them directly, but give an overview of some of the key results needed for our discussion on monoidal monads.  We invite the reader to look in \cite{ReutterZetto} for an excellent summary and further details, but freely use results from other models of enriched categories.

\begin{theorem}[\cite{ReutterZetto} corollary 7.27, corollary 8.20, corollary 8.21]\label{thm:PrtoEnr}
There is a symmetric monoidal functor 
\begin{equation}
	\cat^{(-)}:\alg(\pr)\to \pr
\end{equation}
taking a presentably monoidal category $\mathcal{V}$ to the presentable category $\cat^\mathcal{V}$ of small $\mathcal{V}$-enriched categories. 

\end{theorem}

\begin{example}
	There is an equivalence $\cat\simeq \cat^\mathcal{S}$ \cite[Theorem 5.4.6]{GepnerHaugseng}. 
\end{example}

\noindent It follows from Theorem \ref{thm:PrtoEnr} that a morphism $F:\mathcal{W}\to \mathcal{V}$ in $\alg(\pr)$ induces an adjunction
\begin{equation}
	\begin{tikzcd}
		\cat^\mathcal{\mathcal{W}}\ar[r,shift left=2,"F_!"{name=A}]&\cat^\mathcal{V}\ar[l,shift left=2,"F^*"{name=B}]\ar[phantom, from=A,to=B,"\dashv" rotate=-90]. 
	\end{tikzcd}
\end{equation}
Moreover, 
For $\mathcal{C}\in \cat^ \mathcal{W}$, the category $F_! \mathcal{C}\in \cat^ \mathcal{V}$ has the same objects, and for $x,y\in \mathcal{C}$, the collection of morphisms is given by  
\begin{equation}
	F_!\mathcal{C}(x,y)=F(\mathcal{C}(x,y)). 
\end{equation}

\begin{example}
There is a unique colimit preserving monoidal functor $\mathcal{S}\to \mathcal{V}$, which takes $*\in \mathcal{S}$ to the monoidal unit $\mathds{1}_\mathcal{V}\in \mathcal{V}$.  This induces an adjunction
\begin{equation}
	\begin{tikzcd}
		\cat\simeq \cat^\mathcal{\mathcal{S}}\ar[r,shift left=2,""{name=A}]&\cat^\mathcal{V}\ar[l,shift left=2,"\fun^ \mathcal{V}(\mathds{1}{,}-)"{name=B}]\ar[phantom, from=A,to=B,"\dashv" rotate=-90].
	\end{tikzcd}
\end{equation}
The \emph{underlying category} of $\mathcal{C}\in \cat^ \mathcal{V}$ is its image in $\cat$ under the right adjoint.  We shall denote the underlying category by $\mathcal{C}^{\leq 1}\in \cat$.  Moreover, the underlying category functor is lax monoidal. 
\end{example}

\begin{corollary}\label{corollary:EnrichedMonoidal}
	If $\mathcal{V}\in \alg_{E_1\otimes \mathcal{O}}(\pr)$, then $\cat^\mathcal{V}\in \alg_ \mathcal{O}\pr$ is presentably $\mathcal{O}$-monoidal. 
\end{corollary}

\begin{proof}
	This follows from \ref{thm:PrtoEnr} since monoidal functors take algebras to algebras. 
\end{proof}

\begin{notation}
	Let $\mathcal{V}\in \alg_{E_k}\pr$.  
\begin{itemize}
	\item Denote $0\cat^\mathcal{V}:=\mathcal{V}$.  
	\item For $1\leq n\leq k$, inductively define $n\cat^\mathcal{V}:= \cat^{(n-1)\cat^ \mathcal{V}}$.  This is well defined since $n\cat\in \alg_{E_{k-n}}\pr$ by Corollary \ref{corollary:EnrichedMonoidal}.
	\item $\mathcal{S}$ is symmetric monoidal.  Denote $n\cat:=n\cat^ \mathcal{S}$. 
\end{itemize}
\end{notation}

\begin{definition}\hfill 
	\begin{itemize}
		\item Suppose $k \geq 2$ and $\mathcal{V}\in \alg_{E_k}(\pr)$.  The unit functor $\mathcal{S}\to \mathcal{V}$ induces an adjunction
	\begin{equation}
		\begin{tikzcd}
			n\cat\ar[r,shift left =2,""{name=a}]& n\cat^\mathcal{V}\ar[l,shift left=2, "(-)^{\leq n}"{name=b}]\ar[from=a, to=b, phantom, "\dashv"rotate=-90]
		\end{tikzcd}
	\end{equation}
	for $0\leq n \leq k$.  We may check inductively that the left adjoint is $E_{k-n}$-monoidal.  If $\mathcal{C}\in n\cat^ \mathcal{V}$ is a $\mathcal{V}$-enriched $n$-category, then its \emph{underlying $n$-category} is its image in $n\cat$ under the right adjoint functor above. 
\item  When $0\leq n\leq k-1$, the monoidal functor $n\cat \to n\cat^ \mathcal{V}$ induces an adjunction
	\begin{equation}
		\begin{tikzcd}
			\pr^{n\cat}\ar[r, shift left =2,""{name=a}]&\pr^{n\cat^ \mathcal{V}}\ar[l, shift left =2, "(-)^{\leq n}"{name=b}]\ar[from =a, to=b, phantom, "\dashv"rotate=-90]
		\end{tikzcd}
	\end{equation}
	For $\mathcal{E}\in \pr^{n\cat^ \mathcal{V}}$, we also call its image in $\pr^{n\cat}$ under the above functor the \emph{underlying (large) $n$-category}.  In particular, $n\cat^ \mathcal{V}\in \pr^{n\cat^ \mathcal{V}}$ naturally determines a (large) $n$-category. 
	\end{itemize}
\end{definition}

\begin{example}
	Consider $(n+m)\cat:=n\cat^{m\cat}$.  Then the underlying $n$-category functor $(n+m)\cat\to n\cat$ truncates by discarding all $i$-morphisms with $i>n$.  This motivates the notation $(-)^{\leq n}$. 
\end{example}

\begin{theorem}[\cite{GepnerHaugseng} Theorem 6.32]\label{delooping}
	Let $\mathcal{V}\in \alg_{E_1\otimes \mathcal{O}}(\Pr)$.  The delooping functor $\rB$ is fully faithful and $\mathcal{O}$-monoidal, with right adjoint.  
	\begin{equation}
		\begin{tikzcd}
			\alg(\mathcal{V})\ar[r,shift left =2,"\rB"{name=B}]&\ar[l,shift left=2,"\Omega"{name=A}]\alg_{E_0}(\cat^\mathcal{V})\ar[phantom,from=B,to=A,"\dashv" rotate=-90]
		\end{tikzcd}
	\end{equation}
The right adjoint $\Omega :\cat^ \mathcal{V}_*\to \alg(\mathcal{V})$ is called \emph{looping}, and acts by taking $(\mathcal{C},c)\in \cat^ \mathcal{V}_*$ to the endomorphism object $\mathcal{C}(c,c)$.  
\end{theorem}

\begin{remark}
	Another construction of this adjunction is described in \cite[Remark 13.1.7]{StefanichPhD}.  
\end{remark}

\begin{example}
	Suppose that $\mathcal{V}\in \alg_{E_k}\pr$.  For $1\leq n\leq k$, we may iterate this construction to get an adjunction
\begin{equation}
		\begin{tikzcd}
			\alg_{E_k}(\mathcal{V})\ar[r,shift left =2,"\rB^n"{name=B}]&\ar[l,shift left=2,"\Omega^n"{name=A}]\alg_{E_{k-n}}(n\cat^\mathcal{V})\ar[phantom,from=B,to=A,"\dashv" rotate=-90], 
		\end{tikzcd}
	\end{equation}
	where $\rB^n$ is fully faithful. 
\end{example}

\begin{proposition}[\cite{ReutterZetto}  Proposition 7.3, Observation 8.13,\cite{HeineWeighted} Corollary 3.101]\hfill
	Suppose that $\mathcal{V}$ is presentably $E_1\otimes \mathcal{O}$-monoidal.  The free cocompletion functor $\mathrm{cc}$ is $\mathcal{O}$-monoidal and left adjoint to the inclusion $\Pr^\mathcal{V}\to \widehat{\cat}^\mathcal{V}$.  
\begin{equation}
	\begin{tikzcd}
\widehat{\cat}^\mathcal{V}\ar[r,"\mathrm{cc}"{name=A},shift left=2]&\Pr^\mathcal{V}\ar[l,"\iota"{name=B},shift left=2]\arrow[phantom,from=A,to=B,"\dashv" rotate=-90]
	\end{tikzcd}
\end{equation}
Its restriction to $\cat^\mathcal{V}$ is faithful, and we can identify it with the enriched presheaf functor.  Its unit is the Yoneda embedding, which we denote by $\yo^\mathcal{V}_\mathcal{C}:\mathcal{C}\to \psh^\mathcal{V}(\mathcal{C})$. 
\end{proposition}

\begin{definition}\hfill 
	\begin{itemize}
		\item For $\mathcal{I}\in \cat$, consider the unique morphism $\mathcal{I}\to *$.  Denote the image of this functor under the free functor $\cat\to \cat^ \mathcal{V}$ by $\iota: \mathcal{I}_ \mathcal{V}\to \mathds{1}$. 
		\item Let $F:\mathcal{I}_ \mathcal{V}\to \mathcal{E}$ be a morphism in $\cat^ \mathcal{V}$.  A \emph{conical colimit} for $F$ is a corepresenting object $\colim_ \mathcal{I}F\in \mathcal{E}$ 
			\begin{equation}
				\fun^ \mathcal{V}(\mathcal{I}, \mathcal{E})( \iota^* -,F)\simeq \fun^\mathcal{V}(\mathds{1},\mathcal{E})(\colim_ \mathcal{I}F,-). 
			\end{equation}
	\end{itemize}
	Similarly, a \emph{conical limit} is a representing object $\lim_ \mathcal{I}F\in \mathcal{E}$
	\begin{equation}
		\fun^\mathcal{V}(\mathcal{I},\mathcal{E})(F, \iota^* -)\simeq \fun^ \mathcal{V}(\mathds{1}, \mathcal{E})( -, \lim_ \mathcal{I}F). 
	\end{equation}
\end{definition}

\section{Condensation in enriched categories}\label{condensationsec}
In the following, fix $\mathcal{V}\in \alg_{E_k}\pr$, $n\leq k$, and $\mathfrak{D}\in n\cat^ \mathcal{V}$.  We would like to think of $\mathfrak{D}$ as a \emph{category of $n$-dimensional theories}.  Objects in this category are $n$-dimensional quantum field theories and $1$-morphisms are codimension-$1$ topological interfaces.  More generally, an $i$-morphism is a codimension-$i$ interface between codimension-$(i-1)$ interfaces.  In particular, given a theory $d\in \mathfrak{D}$, the collection of \emph{codimension-$i$ topological defects} in $d$ is given by $\Omega^i_d \mathfrak{D}$.  Before imposing additional properties to ensure that this category behaves well, we shall sketch how one recovers many familiar features from the ordinary fusion categorical theory of condensation in this language.

\begin{example}\label{theorybimodule}
	Consider two objects (`theories') $v,w\in \mathfrak{D}$.  We claim that defects in the theories $v,w$ act on the collection of interfaces between $v$ and $w$.  The collection of codimension-$1$ defects in $v$ is $\Omega_v \mathfrak{D}$, and similarly for $w$.  The collection of interfaces is $\mathfrak{D}(v,w)$.  Restricting the hom-functor gives
	\begin{equation}
		\rB\Omega_v \mathfrak{D}^\op\times \rB \Omega_w \mathfrak{D}\hookrightarrow \mathfrak{D}^\op\times \mathfrak{D}\xrightarrow{\mathfrak{D}(-,-)} (n-1)\cat^\mathcal{V}, 
	\end{equation}
	taking the unique object on the left to the object $\mathfrak{D}(v,w)$.  We get the desired bimodule structure using the equivalence 
\begin{equation}
	\fun(\rB \Omega_v \mathfrak{D}^\op\times \rB \Omega_w \mathfrak{D}, (n-1)\cat^ \mathcal{V})\simeq \bimod_{\Omega_w \mathfrak{D}|\Omega_v \mathfrak{D}}(n-1)\cat^ \mathcal{V}. 
\end{equation}
\end{example}

\begin{definition}\label{defEilenbergMoore}Suppose that $\mathfrak{D}\in 2\cat$ and $y\in \mathfrak{D}$.  Suppose that $T\in \alg(\mathfrak{D}(y,y))$ is a monad on $y\in \mathfrak{D}$.  
	\begin{itemize}
		\item An \emph{Eilenberg-Moore object} for $T$ is a representing object for the functor 
\begin{align}
	\Lmod^*_T(\mathfrak{D}(-,y)):\mathfrak{D}^{1\op}\to \cat&& x\mapsto \Lmod_T( \mathfrak{D}(x,y)).  
		\end{align}
This functor is given precisely in Construction \ref{EMRepresentable}.
\item We say that $\mathfrak{D}$ \emph{admits Eilenberg-Moore objects at $y\in \mathfrak{D}$} if an Eilenberg-Moore object exists for every monad on $y$. 
	\end{itemize}
\end{definition}

\begin{example}\label{EM1}
	  Suppose that $v\in \mathfrak{D}$.  A monad $A\in \alg(\Omega_v \mathfrak{D})$ describes a family of 0-form operators in $v$ which is closed under fusion.  Suppose that $\mathfrak{D}$ admits Eilenberg-Moore objects (definition \ref{defEilenbergMoore}).  Then the left Eilenberg-Moore functor $\EM_L:\alg(\Omega\mathfrak{C})\to \mathfrak{C}_{/v}$ picks out another object $\EM(A)\in \mathfrak{D}$ (another theory!), equipped with an adjunction 
	\begin{equation}
		\begin{tikzcd}
			v\ar[r,shift left=2,"\EM_L(A)"{name=A}]&\EM(A)\ar[l,shift left=2,"\EM_R(A)"{name=B}"]\ar[phantom,from=A,to=B,"\dashv"rotate=-90]. 
		\end{tikzcd}
	\end{equation}
	We will call $\EM(A)$ the \emph{condensate} of $A$.  We now describe the universal property of the Eilenberg-Moore object physically.  If we have $w\in \mathfrak{D}$ and interface $v\leftarrow w$ which is invariant under fusion with $A$, then this interface can be replaced by a thin slice of the condensate $\EM(A)$ in a unique way. 
	\begin{figure}[h]
		\centering
\begin{tikzpicture}
	\draw[thick] (0,0.75)--(0,-.75);
	\draw[thick] (6,.75)--(6,-.75);
	\draw[thick] (8,.75)--(8,-.75);
	\node at (-1,0){$v$};
	\node at (1,0){$w$};
	\node at (0,1.1){$f$};
	\node at (5,0){$v$};
	\node at (7,0){$\EM(A)$};
	\node at (9,0){$w$};
	\node at (6,1.1){$\EM_R(A)$};
	\node at (8,1.1){$\tilde{f}$};
	\node at (3,0){$\simeq$};
\end{tikzpicture}
\end{figure}

The notion of condensation treated by \cite[]{DouglasReutter2018} is given by the Eilenberg-Moore object of a \emph{separable} monad.  We merely treat the case which does not assume separability. 
\end{example}
\begin{remark}
In the example above, while $f\in \mathfrak{D}(w,v)$, $v$ is drawn to the left, so that the left-right ordering of the image is compatible with the ordering of composable morphisms. 
\end{remark}

\subsection{Morita equivalences}\label{MoritaSec}
A typical case in the fusion category literature corresponds to $\mathfrak{D}=\Rmod_A 2\Vect$ being the 2-category of modules over a fusion category $A\in \alg(2\Vect)$.  Objects in this $\Rmod_A 2\Vect$ correspond, informally, to 2-dimensional theories with $A$-symmetry.  The assignment $A\mapsto \Rmod_A 2\Vect$ is an example of the Eilenberg-Moore condensation construction of Example \ref{EM1} (See Corollary \ref{AlgebraEML}).  When everything is sufficiently finite, we expect that all theories can be obtained from a single distinguished theory via condensation.  This is captured for fusion categories by \cite[Theorem 7.10.1]{EGNO}.  

More generally, we could take $\mathfrak{D}=\Rmod_A \mathcal{W}$ for some algebra $A$ in a monoidal category $\mathcal{W}$.  We will now prove Morita-equivalence statements analogous to those proven in \cite[]{EGNO}.  

\begin{notation}
	Let $M\in \Rmod_A \mathcal{W}$.  Denote by $\langle M\rangle \subset \Rmod_A \mathcal{W}$ the full subcategory generated by $M$ under simplicial colimits and tensoring by objects of $\mathcal{W}$. 
\end{notation}
\begin{proposition}\label{prop:MonadMorita}
	Suppose that $\Rmod_A \mathcal{W}$ is a closed $\mathcal{W}$-module category, and $M\in \Rmod_A \mathcal{W}$ is right dualizable and $A\in \langle M\rangle$.  Then 
	\begin{equation}
\Rmod_A\mathcal{W}(M,-)^\mathrm{enh}:\Rmod_A \mathcal{W}\to \Rmod_{\Omega_M \Rmod_A \mathcal{W}} \mathcal{W}
	\end{equation}
	 is an equivalence.  
\end{proposition}

\begin{proof}
	We consider the functor $\Rmod_A \mathcal{W}(M,-): \Rmod_A \mathcal{W}\to \mathcal{W}$ corepresented by $M$.  Our goal is to show that this fulfils the criterion of the $\mathcal{W}$-enriched Barr-Beck-Lurie theorem.  By definition of hom objects, this has a left adjoint, given by $-\otimes_A M$.  By \cite[Proposition 4.6.2.1]{LurieHA}, $\Rmod_A \mathcal{W}(M,-)\simeq -\otimes_A M^\vee$ has a right adjoint given by $\Rmod_A \mathcal{W}(M^\vee, -)$, where $M^\vee$ is the dual to $M$.  Then $\Rmod_A \mathcal{W}(M,-)$ preserves tensors and conical simplicial colimits.  Notice that $\langle A\rangle \simeq\Rmod_A \mathcal{W}$.  Since $A\in \langle M\rangle$, we see that $\langle M\rangle \simeq \Rmod_A \mathcal{W}$.  Since $M$ generates $\Rmod_A \mathcal{W}$ under simplicial colimits and tensors and $\Rmod_A \mathcal{W}(M,-)$ preserves these colimits, we conclude that $\Rmod_A \mathcal{W}(M,-)$ is conservative.  Then $\Rmod_A \mathcal{W}(M,-)$ fulfills the conditions of enriched Barr-Beck-Lurie \cite[Proposition 7.5.9]{StefanichPhD}, we conclude that it is monadic.  The result follows. 
\end{proof}

\begin{remark}
	If $\mathcal{W}$ has enough colimits, then $M\in \Rmod_A \mathcal{W}$ is right dualizable if and only if $\Rmod_A \mathcal{W}(M,-)$ preserves all $\mathcal{W}$-enriched weighted colimits.  That is to say, if $M$ is \emph{atomic} in $\Rmod_A \mathcal{W}$.  See \cite[Lemma 3.22]{ReutterZetto}. 
\end{remark}

\begin{example}
	Let $\mathcal{W}$ be the 2-category of linear categories, and $A\in \alg \mathcal{W}$ be a fusion category.  Then $\Omega_M \Rmod_A \mathcal{W}$ is the (opposite of the) \emph{dual fusion category} $\mathcal{C}^*_M$ with respect to the module $M$ \cite[Definiton 7.12.2]{EGNO}.  Applied to this case, Proposition \ref{prop:MonadMorita} recovers \cite[Theorem 7.12.16]{EGNO}.  

	This gives a Morita equivalence between $A$ and $\Omega_M \Rmod_A \mathcal{W}$.  In particular, it has an inverse: this is analogous to \cite{EGNO} theorem 7.12.11, and recovers the `dual symmetry' of \cite{BhardwajTachikawa2017}.  See also \cite[Secion 3.4.3]{DimofteNiu} for similar physical constructions.  From this perspective, Proposition \ref{prop:MonadMorita} extends the (de)equivariantization or (de)gauging construction of fusion categories.  
\end{example}

\subsection{Condensation of codimension-$i$ operators}\label{sec:iCondense}
We will now impose some additional structure on our `category of theories' $\mathfrak{D}$ in order to get a reasonable theory of condensation of higher-codimension operators.  We will require that they locally have certain colimits, and composition is compatible with these colimits. 

\begin{definition}
	Let $\mathcal{V}\in \alg_{E_k}\pr$, and $\mathcal{K}$ be a collection of categories.  
	\begin{itemize}
		\item Define $\cat^ \mathcal{V}(\coprod^{\geq 0} \mathcal{K})\subset \cat^ \mathcal{V}$ to be the subcategory of $\cat^ \mathcal{V}$ on those enriched categories which admit $\mathcal{K}$-indexed conical colimits, and functors which preserve them.  
		\item For $n\leq k$, inductively denote $n\cat^ \mathcal{V}(\coprod^{\geq 0} \mathcal{K}):= \cat^{(n-1)\cat^ \mathcal{V}(\coprod^{\geq 0} \mathcal{K})}(\coprod \mathcal{K})$.  These are the categories that admit $\mathcal{K}$-indexed colimits, locally admit $\mathcal{K}$-indexed colimits, locally locally, and so on. 
		\item For $n\leq k$, inductively denote $n\cat^ \mathcal{V}(\coprod^{\geq 1} \mathcal{K}):= \cat^{(n-1)\cat^ \mathcal{V}(\coprod^{ \geq 0})}$.  These are the categories which locally admit $\mathcal{K}$-indexed colimits, locally locally ...
	\end{itemize}
\end{definition}

In the following, let $\mathcal{K}$ be a collection of categories containing $\Delta^\op$.  Let $\mathfrak{D}\in n\cat^ \mathcal{V}( \coprod^{\geq 1} \mathcal{K})$ be a category which is locally compatible with $\mathcal{K}$-indexed colimits, and $d\in \mathfrak{D}$.  Given a family of $i$-dimensional operators with an algebra structure satisfying some conditions, we expect to be able to construct a new theory by condensing.  We define this algebraically below.  

\begin{corollary}\label{cor:EMCondensation}
Suppose $A\in \alg_{E_{i+1}}(\Omega^{i+1}_d \mathfrak{D})$ admits an Eilenberg-Moore object.  Then $\EM(A)\in \Omega^i_d \mathfrak{D}$ has a natural $E_i$-algebra structure, which we also denote by $\EM(A)\in \alg_{E_i}(\Omega^i_d \mathfrak{D})$. 
\end{corollary}

\begin{proof}
	This follows immediately from Theorem \ref{EML}.  
\end{proof}

\begin{definition}\hfill 

\begin{itemize}
	\item Let $\alg^c(\Omega \mathfrak{D})\subset \alg(\Omega \mathfrak{D})$ be the full subcategory on those monads which admit Eilenberg-Moore objects in $\mathfrak{D}$.  
	\item For $i\leq n$, let $\alg^c_{E_i}(\Omega^i \mathfrak{D})\subset \alg_{E_i}(\Omega^i \mathfrak{D})$ inductively denote the full subcategory on those algebras $A$ such that 
		\begin{enumerate}
			\item $A$ admits an Eilenberg-Moore objects in $\Omega^{i-1} \mathfrak{D}$. 
			\item The algebra structure $\EM(T)\in \alg_{E_{i-1}} (\Omega^{i-1}\mathfrak{D})$ induced by Corollary \ref{cor:EMCondensation} lands in the subcategory $\alg^c_{E_{i-1}}(\Omega^{i-1} \mathfrak{D})$. 
		\end{enumerate}
	\item For $l\leq i\leq n$ and $A\in \alg^c_{E_i}( \mathfrak{D})$, inductively denote $\EM^l(A):= \EM(\EM^{l-1}(A))$.  When $l<i$, this is an object in $\alg^c_{E_{i-l}}(\Omega^{l-i} \mathfrak{D})$.  When $l=i$, it is an object in $\mathfrak{D}$. 
\end{itemize}
The algebras $T\in \alg^c_{E_i}(\Omega^i \mathfrak{D})$ will be called \emph{condensable}.  It is immediate from this construction that, given $d\in \mathfrak{D}$ and $A\in \alg^c_{E_i}(\Omega^i \mathfrak{D})$, one can construct a new pointed object $\EM^i(A)\in \mathfrak{D}$ called the \emph{condensate}.  
\end{definition}

\subsection{Centers, centralizers, and $E_k$-modules}\label{centralizersec}

Below, we will introduce the notions of center and centralizer for \emph{coherent} \cite[Definition 3.3.1.9]{LurieHA} operads.  We will primarily care for the case $\mathcal{O}=E_k$, which does indeed satisfy these conditions. 
\begin{definition}[\cite{LurieHA} 5.3.1.2, 5.3.1.12]
Suppose that $\mathcal{O}$ is a coherent symmetric $\infty$-operad and $\mathcal{W}\in \alg_{E_1\otimes \mathcal{O}}(\pr)$.  Then $\alg_\mathcal{O}(\mathcal{W})$ is $E_1$-monoidal, hence a left module over itself.  Let $\mathds{1}_\mathcal{W}\in \alg_\mathcal{O}(\mathcal{W})$ be the monoidal unit, and suppose we have morphism $f:\mathds{1}_\mathcal{W}\otimes C\to B$ in $\alg_\mathcal{O}(\mathcal{W})$.  The \emph{$\mathcal{O}$-centralizer} of $f$ is terminal amongst objects $A\in \alg_\mathcal{O}(\mathcal{W})$ equipped with a morphism $A\otimes C\to B$ making the following diagram commute
	\begin{equation}
		\begin{tikzcd}
			& A\otimes C\ar[dr]&\\
			\mathds{1}_\mathcal{W}\otimes C\ar[ur,"\iota_A\otimes \id_C"]\ar[rr,swap,"f"]&&B,
		\end{tikzcd}
	\end{equation}
	where $\iota_A:\mathds{1}_\mathcal{W}\to A$ is the unit morphism.  The $\mathcal{O}$-centralizer of $f$ exists by \cite{LurieHA} theorem 5.3.1.14, and will be denoted by $\FZ_\mathcal{O}(f)$. 	
\end{definition}

\begin{definition}
	Suppose that $\mathcal{W}\in \alg_{E_1\otimes \mathcal{O}}(\pr)$ and $C\in \alg_\mathcal{O}(\mathcal{W})$.  The \emph{$\mathcal{O}$-center} of $C$ is the $\mathcal{O}$ centralizer of $\id_C:C\to C$, denoted by $\FZ_\mathcal{O}(C):=\FZ_\mathcal{O}(\id_C)\in \alg_{E_1\otimes \mathcal{O}}(\mathcal{W})$, where we note that $\FZ_\mathcal{O}(\id_C)$ obtains an additional $E_1$-monoidal structure as in \cite[Remark 5.3.1.13]{LurieHA}. 
\end{definition}
\begin{example}
	Consider $\cat\in \alg_{E_\infty}(\pr)$. 
	\begin{itemize}
		\item Consider a monoidal category $\mathcal{C}\in \alg(\cat)$.  The \emph{Drinfeld center} of $\mathcal{C}$ is $\FZ_{E_1}(\mathcal{C})\in \alg_{E_2}(\cat)$.  
		\item Consider a braided monoidal category $\mathcal{D}\in \alg_{E_2}(\cat)$.  The \emph{M\"uger center} of $\mathcal{D}$ is $\FZ_{E_2}(\mathcal{D})\in \alg_{E_3}(\cat)$. 
	\end{itemize}
	
\end{example}
\begin{theorem}[\cite{LurieHA} 5.3.1.30]\label{5.3.1.30}
	Suppose that $\mathcal{O}$ is a coherent symmetric $\infty$-operad, and $\mathcal{W}\in \alg_{E_1\otimes \mathcal{O}}(\pr)$ and $f:C\to D$ is a morphism in $\alg_\mathcal{O}(\mathcal{W})$.  Define $\Mod^{\mathcal{O}}_C(\mathcal{W})$ as in \cite{LurieHA} 3.3.3, and consider $f^*D\in \Mod^{\mathcal{O}}_{\mathcal{C}}(\mathcal{W})$ (\cite{LurieHA} 3.4.1.7).  Then we have an equivalence 
	\begin{equation}
		\FZ_\mathcal{O}(f)\simeq \Mod^{\mathcal{O}}_C(\mathcal{W})(C,f^*D)
	\end{equation}
	between the $\mathcal{O}$-centralizer of $f$ and the morphism-object between $C$ and $D$, viewed as $\mathcal{O}$-modules over $C$. 
\end{theorem}
\begin{corollary}
The $\mathcal{O}$-center of $C$ can be identified with endomorphisms of $C$, viewed as an $\mathcal{O}$-module over itself. 
\begin{equation}
	\FZ_\mathcal{O}(C)\simeq \Omega_ C\Mod^\mathcal{O}_C(\mathcal{W})
\end{equation}

\end{corollary}
The following is adapted from the proof of \cite{LurieHA} 5.3.2.5.  
\begin{proposition}\label{iterativecentralizer}

	Suppose that $\mathcal{W}\in \alg_{E_{k+1}}(\pr)$, and $f:C\to D$ is a morphism in $\alg_{E_k}(\mathcal{W})$.  Then one has an equivalence in $\alg_{E_k}(\mathcal{W})$
\begin{equation}
	\FZ_{E_k}(f)\simeq \Omega \FZ_{E_{k-1}}(\Rmod_f \mathcal{W}). 
\end{equation}
\end{proposition}

\begin{proof*}
Using corollary \ref{AlgebraEML} (or alternatively \cite{LurieHA} theorem 4.8.5.16) applied to a larger universe, together with \cite{LurieHA} 4.8.5.11, one obtains an adjunction 
\begin{equation}
\begin{tikzcd}
	Alg_{E_k}(\mathcal{W})\ar[r,shift left=2,"\Rmod_{-} \mathcal{W}"]&Alg_{E_{k-1}}(\pr^\mathcal{W})\ar[l,shift left=2,"\Omega"]
\end{tikzcd}
\end{equation}
where $\Rmod_{!-}(\mathcal{V})$ is fully faithful and monoidal.  For any $A\in \alg_{E_k}(\mathcal{W})$, this gives equivalence 

\begin{equation}
	\alg_{E_k}(\mathcal{W})(A,\Omega \FZ_{E_{k-1}}(f))\simeq \alg_{E_{k-1}}(\pr^\mathcal{W})(\Rmod_A \mathcal{W},\FZ_{E_{k-1}}(f)). 
\end{equation}
Using the universal property of $\FZ_{E_{k-1}}(\Rmod_f \mathcal{W})$, this is identified with the collection of morphisms $\Rmod_A \mathcal{W}\otimes \Rmod_C \mathcal{W}\to \Rmod_D \mathcal{W}$ making the following diagram commute 
\begin{equation}
		\begin{tikzcd}
			&\Rmod_A \mathcal{W}\otimes \Rmod_C \mathcal{W}\ar[dr]&\\
			\Rmod_C \mathcal{W}\ar[ur]\ar[rr,swap,"\Rmod_f \mathcal{W}"]&&\Rmod_ D \mathcal{W}. 
		\end{tikzcd}
	\end{equation}
	Using that $\Rmod_{-} \mathcal{W}$ is fully faithful and monoidal, this is identified with the collection of morphisms $A\otimes C\to D$ making the following diagram commute

\begin{equation}
		\begin{tikzcd}
			&A\otimes C\ar[dr]&\\
			C\ar[ur]\ar[rr,swap,"f"]&&D, 
		\end{tikzcd}
	\end{equation}
	which, using the universal property of $\FZ_{E_k}(f)$, is in turn identified with $\alg_{E_k}(\mathcal{W})(A,\FZ_{E_k}(f))$.  We then get a natural isomorphism $\alg_{E_k}(\mathcal{W})(A,\FZ_{E_k}(f))\simeq \alg_{E_{k-1}}(\pr^\mathcal{W})(A,\Omega \FZ_{E_{k-1}}(\Rmod_f \mathcal{W}))$.  The result follows from Yoneda lemma. 
\end{proof*}
This offers a simplified method of calculating centralizers, since $E_1$ and $E_0$ centralizers are more easily accessible using Theorem \ref{5.3.1.30}.  Up to some set-theoretic issues which are addressed in \cite[Section 12.2]{StefanichPhD} (see also \cite[Section 8.4]{StefanichPhD} and \cite[]{Aoki}), this construction can be iterated.  For $A\in \alg_{E_k}(\mathcal{W})$ and $i\leq k$, we informally inductively define $\Rmod^{i+1}_A:= \Rmod_{\Rmod^i_A}$.  We recover \cite{Kong2024} theorem 2.3.25. 

\begin{corollary}\label{cor:IterativeCenter}
	Suppose that $\mathcal{W}\in \alg_{E_{k+1}}(\pr)$ and $C\in \alg_{E_k}(\mathcal{W})$.  Then one has an equivalence in $\alg_{E_{k+1}}(\mathcal{W})$
	\begin{equation}
		\FZ_{E_k}(C)\simeq \Omega^{k-i}\FZ_{E_i}(\Rmod_C^{k-i}\mathcal{W} ). 
	\end{equation}
	
\end{corollary}

\begin{conjecture}\label{conj:EiMod}
	There is an equivalence 
	\begin{equation}
		\Omega\Mod^{E_i}_{\Rmod_A}\simeq \Mod^{E_{i+1}}_A
	\end{equation}
\end{conjecture}

\begin{corollary}\label{cor:HigherEiMod}
	Suppose that $A\in \alg_{E_{i+n}}$.  There is an equivalence 
\begin{equation}
		\Omega^n\Mod^{E_i}_{\Rmod^n_A}\simeq \Mod^{E_{i+n}}_A
	\end{equation}
\end{corollary}

\begin{proof}
	We prove by induction on $n$.  The case $n=1$ is Conjecture \ref{conj:EiMod}.  Then 
	\begin{equation}
		\Omega^{n}\Mod^{E_i}_{\Rmod^{n}_A}\simeq \Omega^{n-1}\Omega \Mod^{E_i}_{\Rmod_{\Rmod^{n-1}_A}}\simeq \Omega^{n-1}\Mod^{E_{i+1}}_{\Rmod^{n-1}_A}. 
	\end{equation}
	The result follows from the inductive hypothesis. 
\end{proof}

\begin{corollary}\label{cor:IterativeEn}
	Suppose that $A\in \alg_{E_{i+n}}$.  There is an equivalence 
	\begin{equation}
		\Omega^{n+1}_{\Rmod^{n+1}_A}\simeq \Mod^{E_{n+1}}_A. 
	\end{equation}
\end{corollary}

\begin{proof}
	This follows from Corollary \ref{cor:HigherEiMod}, together with \cite[Remark 4.8.4.9]{LurieHA}, noting that $\bimod_{A|A}\simeq \Mod^{E_1}_A$. 

\begin{equation}
	\Omega^{n+1}_{\Rmod^{n+1}_A}\simeq \Omega^n\left(\Omega_{\Rmod_{\Rmod^n_A}}\right)\simeq \Omega^n \Mod^{E_1}_{\Rmod^n_A}\simeq \Mod^{E_{n+1}}_A. 
	\end{equation}
\end{proof}

\noindent This also gives an alternate proof of Corollary \ref{cor:IterativeCenter}. 
\begin{corollary}
	Suppose that $A\in \alg_{E_{i+n}}$.  There is an equivalence $\FZ_{E_{n+1}}(A)\simeq \Omega \FZ_{E_n}(\Rmod_A)$. 
\end{corollary}

\begin{proof}
	This follows from Corollary \ref{cor:IterativeEn} together with \cite[Theorem 5.3.1.30]{LurieHA}. 
	\begin{equation}
		\FZ_{E_{n+1}}(A)\simeq \Omega \Mod^{E_{n+1}}_A\simeq \Omega^{n+2}_{\Rmod^{n+1}_A}\simeq \Omega \Omega^{n+1}_{\Rmod^n_{\Rmod_A}}\simeq \Omega \FZ_{E_n}(\Rmod_A) 
	\end{equation}

\end{proof}

\section{Construction of the monoidal Eilenberg-Moore functors}\label{Monadssec}
In the last section, we will construct the monoidal Eilenberg-Moore functors that were used to define condensation of codimension-$i$ operaters in Section \ref{sec:iCondense}.  This will require introducing some additional mathematical machinery.

\subsection{Day convolution}\label{Dayconvolution}
In the following, we prove functoriality of the enriched Day convolution, introduced in \cite{Hinich2023} section 8.3, by generalizing the approach of \cite{BenMosheSchlank} to the enriched setting.  We refer to \cite[]{BlomLoubatonRuit,Winges2025} and \cite[Section 2.2.6]{LurieHA} for the unenriched setting. 

\begin{definition}
	Suppose that $\mathcal{V}\in \alg_{E_1\otimes \mathcal{O}}(\pr)$ so that $\cat^ \mathcal{V}\in \alg_ \mathcal{O}(\pr)$.  An \emph{$\mathcal{O}$-monoidal $\mathcal{V}$-category} is an object $\mathcal{C}\in \alg_ \mathcal{O}(\cat^ \mathcal{V})$.  Notice that the underlying category functor $\cat^ \mathcal{V}\to \cat$ is lax $\mathcal{O}$-monoidal, so the underlying category of an $\mathcal{O}$-monoidal $\mathcal{V}$-category is also $\mathcal{O}$-monoidal. 
\end{definition}

\begin{theorem}\label{thm:DayConv}\cite[Section 8.3, 8.4.2]{Hinich2023}
Suppose that $\mathcal{V}\in \alg_{E_2\otimes \mathcal{O}}(\pr)$, $\mathcal{B}\in \alg_\mathcal{O}(\pr^\mathcal{V})$, and $\mathcal{C},\mathcal{D}\in \alg_\mathcal{O}(\cat^\mathcal{V})$.  There is an object $\fun^\mathcal{V}(\mathcal{C},\mathcal{B})^\otimes\in \alg_ \mathcal{O}(\pr^\mathcal{V})$ whose underlying presentably enriched category is $\fun^\mathcal{V}(\mathcal{C},\mathcal{B})$.  Moreover, there is an $\mathcal{O}$-monoidal $\mathcal{V}$-linear equivalence  
	\begin{equation}
		\fun^\mathcal{V}(\mathcal{D},\fun^\mathcal{V}(\mathcal{C},\mathcal{B})^\otimes)^\otimes \simeq \fun^\mathcal{V}(\mathcal{D}\otimes \mathcal{C},\mathcal{B})^\otimes. 
	\end{equation}
\end{theorem}

\begin{definition}
 Suppose that $\mathcal{C}\in \alg_ \mathcal{C}\cat^\mathcal{V}$ and $\mathcal{B}\in \alg_ \mathcal{O}\pr^\mathcal{V}$.  
	\begin{itemize}
		\item An \emph{$\mathcal{O}$-algebra object} in $\mathcal{C}$ is an $\mathcal{O}$-algebra object in the underlying monoidal category of $\mathcal{C}$. 
		\item The monoidal structure $\fun^\mathcal{V}(\mathcal{C},\mathcal{B})^\otimes$ described in Theorem \ref{thm:DayConv} will be called the \emph{Day convolution} monoidal structure on $\fun^\mathcal{V}(\mathcal{C},\mathcal{B})$.  
		\item A \emph{lax monoidal functor} from $\mathcal{C}$ to $\mathcal{B}$ is an object of $\fun^\mathcal{V}(\mathcal{C},\mathcal{B})_{ \lax \mathcal{O}}:= \alg_ \mathcal{O}\fun^\mathcal{V}(\mathcal{C},\mathcal{B})^\otimes$.  
	\end{itemize}
\end{definition}
\begin{remark}\label{rmk:LaxMonModFunc}
	In the case $\mathcal{V}=\mathcal{S}$, it follows from any of \cite[]{BlomLoubatonRuit,Winges2025} or \cite[Section 2.2.6]{LurieHA} that there is an equivalence 
	\begin{equation}
		\fun(\mathcal{C},\mathcal{B})_{\mathcal{O}\lax}\simeq \opd_{/\mathcal{O}}(\mathcal{C},\mathcal{B}), 
	\end{equation}
	where $\mathcal{C},\mathcal{B}$ are viewed as operads over $\mathcal{O}$ by unstraightening.

	In particular, take $\mathcal{O}=\mathrm{LM}\otimes \mathcal{O}'$, where $\mathrm{LM}$ is the left module operad (\cite[Definition 4.2.1.1]{LurieHA}).  An $\mathcal{O}'$-monoidal presentable $\mathcal{V}$-module is $\mathcal{B}\in \alg_{\mathcal{O}'}(\pr^\mathcal{V})$.  Then $\mathcal{B}$ determines an $\mathcal{O}'\otimes \mathrm{LM}$-algebra in $\pr$.  For any two such presentably monoidal module categories, one can then speak of morphisms of operads over $\mathrm{LM}\otimes \mathcal{O}'$ between them which are strong $\mathrm{LM}$-monoidal, and this agrees with the notion of lax monoidal functor described above.
\end{remark}

\begin{lemma}\label{presentablepostcomp}
	Suppose that that $\mathcal{C},\mathcal{D},\mathcal{E}\in \alg_ \mathcal{O}\cat^\mathcal{V}$ are such that the Day convolutions $\fun_ \mathcal{V}(\mathcal{C},\mathcal{D})^\otimes$ and $\fun_ \mathcal{V}(\mathcal{C},\mathcal{D})^\otimes$ exist, and $g:\mathcal{D}\to \mathcal{E}$ is an $\mathcal{O}$-monoidal functor.  Then the postcomposition functor 
	\begin{equation}
		g\circ -:\fun^ \mathcal{V}(\mathcal{C},\mathcal{D})^\otimes \to \fun^\mathcal{V}(\mathcal{C},\mathcal{E})^\otimes
	\end{equation}
is $\mathcal{O}$-monoidal. 
\end{lemma}

\begin{proof}
	By the monoidal tensor-hom adjunction, $\fun^\mathcal{V}(\mathcal{C},\mathcal{D})^\otimes$ is a representing object for the lax $\mathcal{O}$-monoidal functor $\fun^\mathcal{V}(-\otimes \mathcal{C},\mathcal{B})^\otimes:\cat^\mathcal{V}\to \fun((\cat^\mathcal{V})^\op, \cat^\mathcal{V})$.  Since $g:\mathcal{D}\to \mathcal{E}$ induces an $\mathcal{O}$-monoidal natural transformation, we see that it also induces a morphism $g\circ -:\fun^\mathcal{V}(\mathcal{C},\mathcal{D})^\otimes \to \fun^\mathcal{V}(\mathcal{C},\mathcal{E})^\otimes$ between representing objects. 

\end{proof}

\begin{remark}
The equivalence $\fun^\mathcal{V}(\mathcal{C}^\op,\mathcal{V})\simeq \psh^\mathcal{V}(\mathcal{C})$ lifts to an $\mathcal{O}$-monoidal $\mathcal{V}$-linear equivalence between the Day convolution monoidal structure and the monoidal structure induced by monoidality of the presheaf functor. 
\end{remark}

\begin{lemma}\label{YonedaLeftKan}
	Suppose that $\mathcal{V}\in \alg_{E_2\otimes \mathcal{O}}(\pr)$ and $f:\mathcal{C}\to \mathcal{D}$ is a morphism in $\alg_\mathcal{O}(\cat^\mathcal{V})$.  Then the left Kan extension functor 
	\begin{equation}
		f_!:\fun^\mathcal{V}(\mathcal{C},\mathcal{V})^\otimes\to \fun^\mathcal{V}(\mathcal{D},\mathcal{V})^\otimes
	\end{equation}
is a morphism in $\alg_\mathcal{O}(\pr^\mathcal{V})$, i.e., a strong $\mathcal{O}$-monoidal colimit preserving $\mathcal{V}$-functor. 
\end{lemma}

\begin{proof}
	Under the identification $\fun^\mathcal{V}(\mathcal{C}^{\op},\mathcal{V})^\otimes\simeq \psh^\mathcal{V}(\mathcal{C})$, one can identify $g_!\simeq \psh^\mathcal{V}(g):\psh^\mathcal{V}(\mathcal{C})\to \psh^\mathcal{V}(\mathcal{D})$, which is an $\mathcal{O}$-algebra morphism in $\pr^\mathcal{V}$ by $\mathcal{O}$-monoidality of $\psh^\mathcal{V}:\cat^\mathcal{V}\to \pr^\mathcal{V}$.  We apply this to $f^{\op}:\mathcal{C}^{\op}\to \mathcal{D}^{\op}$. 
\end{proof}

\begin{lemma}\label{YoyoLeftKan}
	Suppose that $\mathcal{V}\in \alg_{E_2\otimes \mathcal{O}}(\pr)$, $\mathcal{C},\mathcal{D},\mathcal{E}\in \alg_\mathcal{O}(\cat^\mathcal{V})$, and $f:\mathcal{D}\to \mathcal{E}$ is a morphism in $\alg_\mathcal{O}(\cat^\mathcal{V})$.  Then the left Kan extension functor
	\begin{equation}
		f_!:\fun^\mathcal{V}(\mathcal{E},\fun^\mathcal{V}(\mathcal{C},\mathcal{V})^\otimes)^\otimes \to \fun^\mathcal{V}(\mathcal{D},\fun^\mathcal{V}(\mathcal{C},\mathcal{V})^\otimes)^\otimes
	\end{equation}
	is $\mathcal{O}$-monoidal. 
\end{lemma}

\begin{proof}
$f_!$ is equivalent to the composition
\begin{equation}
	\fun^\mathcal{V}(\mathcal{D},\fun^\mathcal{V}(\mathcal{C},\mathcal{V})^\otimes)^\otimes\simeq \fun^\mathcal{V}(\mathcal{D}\otimes \mathcal{C},\mathcal{V})
	\xrightarrow{(f\otimes id_\mathcal{C})_!}\fun^\mathcal{V}(\mathcal{E}\otimes \mathcal{C},\mathcal{V})^\otimes
	\simeq \fun^\mathcal{V}(\mathcal{E},\fun^\mathcal{V}(\mathcal{C},\mathcal{V})^\otimes)^\otimes. 
\end{equation}
The result follows from lemma \ref{YonedaLeftKan}. 

\end{proof}

\begin{proposition}\label{reflectivelocalization}
	Suppose that $\mathcal{V}\in \alg_{E_2\otimes \mathcal{O}}(\pr)$ and $\mathcal{B}\in \alg_\mathcal{O}(\pr^\mathcal{V})$.  Then $\mathcal{B}$ is an $\mathcal{O}$-monoidal $\mathcal{V}$-linear localization of $\fun^\mathcal{V}(\mathcal{C},\mathcal{V})^\otimes$ for some $\mathcal{C}\in \alg_\mathcal{O}(\cat^\mathcal{V})$. 
\end{proposition}
\begin{proof}
	Note that $\mathcal{B}$ is a presentably $\mathrm{LM}\otimes\mathcal{O}$-monoidal category.  It follows from \cite{NikolausSagave} theorem 2.2 that we can find small $\mathrm{LM}\otimes\mathcal{O}$-monoidal category $\mathcal{M}$ and $\mathrm{LM}\otimes \mathcal{O}$-monoidal localization $\fun(\mathcal{M},\mathcal{S})^\otimes\to \mathcal{B}^\otimes$.  Looking at the fiber over $E_1\otimes \mathcal{O}\subset \mathrm{LM}\otimes \mathcal{O}$, we find an $E_1\otimes \mathcal{O}$-monoidal category $\mathcal{A}$ with $E_1\otimes \mathcal{O}$-monoidal localization $\fun(\mathcal{A},\mathcal{S})^\otimes \to \mathcal{V}^\otimes$.  Extension of scalars along this localization gives a presentably $\mathcal{O}$-monoidal left $\mathcal{V}$-module structure on $\fun(\mathcal{M},\mathcal{S})\otimes \mathcal{V}$, which we identify with $\fun^\mathcal{V}(\mathcal{M},\mathcal{V})$ using \cite[Lemma 2.25]{ReutterZetto}, where $\mathcal{M}$ has the free $\mathcal{V}$-enrichment.  The $\mathrm{LM}\otimes \mathcal{O}$-monoidal localization $\fun^\mathcal{V}(\mathcal{M},\mathcal{V})\to \mathcal{B}$ is as desired. 
\end{proof}

\begin{proposition}\label{prop:LKEDayConvMon}
	Suppose that $\mathcal{V}\in \alg_{E_2\otimes \mathcal{O}}(\pr)$, $\mathcal{B}\in \alg_\mathcal{O}(\pr^\mathcal{V})$, and $f:\mathcal{C}\to \mathcal{D}$ is a morphism in $\alg_\mathcal{O}(\cat^\mathcal{V})$.  Then the left Kan extension functor  
	\begin{equation}
		f_!:\fun^\mathcal{V}(\mathcal{C},\mathcal{B})^\otimes \to \fun^\mathcal{V}(\mathcal{D},\mathcal{B})^\otimes
	\end{equation}
is $\mathcal{O}$-monoidal. 
\end{proposition}
\begin{proof}
	The proof is identical to \cite{BenMosheSchlank} Proposition 3.6.  By Proposition \ref{reflectivelocalization}, $\mathcal{B}$ is an $\mathcal{O}$-monoidal reflective localization of $\fun^\mathcal{V}(\mathcal{E},\mathcal{B})^\otimes$ for some $\mathcal{E}\in \alg_\mathcal{O}(\cat^\mathcal{V})$.  Then then the functor $$\fun^\mathcal{V}(\mathcal{C},\fun^\mathcal{V}(\mathcal{E},\mathcal{V})^\otimes)^\otimes\to \fun^\mathcal{V}(\mathcal{C},\mathcal{B})^\otimes$$ given by postcomposition is a $\mathcal{V}$-linear reflective localization, and $\mathcal{O}$-monoidal by lemma \ref{presentablepostcomp}, and similarly for $\fun^\mathcal{V}(\mathcal{D},\fun^\mathcal{V}(\mathcal{E},\mathcal{V})^\otimes)^\otimes\to \fun^\mathcal{V}(\mathcal{D},\mathcal{B})^\otimes$.  By lemma \ref{YoyoLeftKan}, we get $\mathcal{O}$-monoidal $\mathcal{V}$-linear functor
	\begin{equation}
		\fun^\mathcal{V}(\mathcal{C}, \fun^\mathcal{V}(\mathcal{E},\mathcal{V})^\otimes)^\otimes \xrightarrow{f_!}\fun^\mathcal{V}(\mathcal{D}, \fun^\mathcal{V}(\mathcal{E},\mathcal{V})^\otimes)^\otimes\to \fun^\mathcal{V}(\mathcal{D},\mathcal{B})^\otimes, 
	\end{equation}
	which factors $\mathcal{O}$-monoidally through the $\mathcal{O}$-monoidal localization $\fun^\mathcal{V}(\mathcal{C},\fun^\mathcal{V}(\mathcal{E},\mathcal{V})^\otimes)^\otimes \to \fun^\mathcal{V}(\mathcal{C},\mathcal{B})^\otimes$, giving the desired $\mathcal{O}$-monoidal functor $f_!:\fun^\mathcal{V}(\mathcal{C},\mathcal{B})^\otimes\to \fun^\mathcal{V}(\mathcal{D},\mathcal{B})^\otimes$. 
\end{proof}

\begin{corollary}\label{precompmon}
	Suppose that $\mathcal{V}\in \alg_{E_2\otimes \mathcal{O}}(\pr)$, $\mathcal{B}\in \alg_\mathcal{O}(\pr^\mathcal{V})$, and $f:\mathcal{C}\to \mathcal{D}$ is a morphism in $\alg_\mathcal{O}(\cat^\mathcal{V})$.  Then the restriction functor 
	\begin{equation}
		-\circ f:\fun^\mathcal{V}(\mathcal{D},\mathcal{B})^\otimes \to\fun^\mathcal{V}(\mathcal{C},\mathcal{B})^\otimes
	\end{equation}
is lax $\mathcal{O}$-monoidal. 	
\end{corollary}

\begin{proof}
	By Proposition \ref{prop:LKEDayConvMon} and Remark \ref{rmk:LaxMonModFunc} and the fact that the right adjoint to a functor of monoidal categories is automatically lax monoidal, we see that there is a lax $\mathcal{O}$-monoidal lax $\mathcal{V}$-linear enhancement of the right adjoint of underlying categories.  By assumption, this functor is $\mathcal{V}$-linear.  The result follows. 
\end{proof}

\begin{proposition}
	$\mathcal{V}\in \alg_{E_2\otimes \mathcal{O}}(\pr)$, $\mathcal{B}\in \alg_\mathcal{O}(\pr^\mathcal{V})$, and $\mathcal{C}\in \alg_\mathcal{O}(\cat^\mathcal{V})$.  Suppose that $f:\mathcal{C}\to \mathcal{B}$ is a lax $\mathcal{O}$-monoidal $\mathcal{V}$-functor.  Then the underlying functor $f^{\leq 1}:\mathcal{C}^{\leq 1}\to \mathcal{B}^{\leq1}$ is lax $\mathcal{O}$-monoidal. 
\end{proposition}
\begin{proof}
	By applying the underlying category functor to the evaluation map $\fun^\mathcal{V}(\mathcal{C},\mathcal{B})^\otimes \mathcal{C}\to \mathcal{B}$, we obtain an evaluation map $\fun^\mathcal{V}(\mathcal{C},\mathcal{B})^{\otimes \leq 1}\otimes \mathcal{C}^{\leq 1}\to \mathcal{B}^{\leq 1}$.  The universal property of $\mathcal{S}$-enriched Day convolution induces a $\mathcal{O}$-monoidal functor $\fun^\mathcal{V}(\mathcal{C},\mathcal{B})^{\otimes\leq 1}\to \fun(\mathcal{C}^{\leq 1},\mathcal{B}^{\leq 1})^\otimes.$  In particular, $f\in \alg_\mathcal{O}(\fun^\mathcal{V}(\mathcal{C},\mathcal{B}))^\otimes:=\alg_\mathcal{O}(\fun_\mathcal{V}(\mathcal{C},\mathcal{B})^\otimes)$ is taken to an $\mathcal{O}$-algebra in $\fun(\mathcal{C},\mathcal{B})^\otimes$, determining a lax $\mathcal{O}$-monoidal functor by Remark \ref{rmk:LaxMonModFunc}.
\end{proof}

\begin{corollary}
	Suppose that $\mathcal{V}\in \alg_{E_2\otimes \mathcal{O}}(\pr)$, $\mathcal{B}\in \alg_\mathcal{O}(\pr^\mathcal{V})$, $\mathcal{C}\in \alg_\mathcal{O}(\cat^\mathcal{V})$, and $f:\mathcal{C}\to \mathcal{B}$ is a lax $\mathcal{O}$-monoidal $\mathcal{V}$-functor.  Then $f$ takes $\mathcal{O}$-algebra objects in $\mathcal{C}$ to $\mathcal{O}$-algebra objects in $\mathcal{B}$.  
\end{corollary}

\begin{corollary}
	Suppose that $\mathcal{V}\in \alg_{E_2\otimes \mathcal{O}}(\pr)$, $\mathcal{B}\in \alg_{\mathcal{O}}(\pr^\mathcal{V})$, $\mathcal{C},\mathcal{D}\in \alg_\mathcal{O}(\cat^\mathcal{V})$, $f:\mathcal{C}\to \mathcal{D}$ is a morphism in $\alg_\mathcal{O}(\cat^\mathcal{V})$, and $g:\mathcal{D}\to \mathcal{B}$ is a lax $\mathcal{O}$-monoidal $\mathcal{V}$-functor.  Then $g\circ f:\mathcal{D}\to \mathcal{B}$ is a lax $\mathcal{O}$-monoidal $\mathcal{V}$-functor. 
\end{corollary}

\begin{proof}
	$-\circ f:\fun^\mathcal{V}(\mathcal{D},\mathcal{B})^\otimes \to \fun^\mathcal{V}(\mathcal{C},\mathcal{B})^\otimes$ is a lax $\mathcal{O}$-monoidal functor of presentable $\mathcal{V}$-modules, and so takes $g\in \alg_\mathcal{O}(\fun^\mathcal{V}(\mathcal{D},\mathcal{B})^\otimes)$ to $g\circ f\in \alg_\mathcal{O}(\fun^\mathcal{V}(\mathcal{C},\mathcal{B})^\otimes)$. 
\end{proof}

\subsection{The right Eilenberg-Moore functor}
Monads determined by adjunctions have some additional structure: they act on the corresponding left adjoint morphism.  First, we demonstrate how the endomorphism category of $y$ acts on categories of morphisms out of $y$.  We also show that this construction respects additional monoidal structure.  The following treatment of monads loosely follows that given by \cite{StefanichPhD} 7.1.3, 7.2.4, 7.3.7.  See also \cite{HeineMonadicity} for another approach.

\begin{construction}\label{Yonedamodule}
	Suppose that $\mathcal{O}$ is a symmetric operad, $\mathcal{V}\in \alg_{E_2\otimes\mathcal{O}}(\pr)$, $\mathcal{D}\in \alg_\mathcal{O}(\cat^\mathcal{V})$ and $y\in \alg_\mathcal{O}\mathcal{D}$.  By delooping hypothesis (Theorem \ref{delooping}), $\Omega \mathcal{D}\in \alg_{E_1\otimes\mathcal{O}}(\mathcal{V})$.  The restriction functor $\fun^\mathcal{V}(\mathcal{D}^\op,\mathcal{V})\to \fun^\mathcal{V}(\rB\Omega \mathcal{D}^\op,\mathcal{V})$ along the inclusion is lax $\mathcal{O}$-monoidal $\mathcal{V}$-linear by corollary \ref{precompmon}.  With monoidal Yoneda embedding, we then have composition of lax $\mathcal{O}$-monoidal $\mathcal{V}$-functors
\begin{equation}
	\mathcal{D}^{enh}(-,y):\mathcal{D}^{1\op}\to \fun^\mathcal{V}(\mathcal{D},\mathcal{V})\to \fun^\mathcal{V}(\rB\Omega \mathcal{D},\mathcal{V})\simeq \Lmod_{\Omega\mathcal{D}}(\mathcal{V}),
\end{equation}
taking $x\in \mathcal{D}$ to $\mathcal{D}(x,\mathds{1})$, viewed as a left $\Omega \mathcal{D}:= \mathcal{D}(\mathds{1},\mathds{1})$-module.  
\end{construction}

\begin{definition}
	Suppose that $f:x\to y$ is a morphism in $\mathfrak{D}\in 2\cat$.  Consider $\mathfrak{D}(x,y)$ as a left module over $\mathfrak{D}(y,y)$.  Then an \emph{endomorphism monad} for $f$ is an endomorphism object $\End(f)\in \mathfrak{D}(y,y)$ for $f$.  
\end{definition}

\begin{proposition}(\cite{StefanichPhD} Proposition 7.2.4)\label{adjendmon}
	Suppose that $\mathfrak{D}\in 2\cat$ and $r:x\to y$ is a morphism in $\mathfrak{D}$.  Suppose that $r$ admits a left adjoint $l:y\to x$.  Then $r$ admits an endormorphism monad with underlying object $rl\in \mathfrak{D}(y,y)$, and action map $rlr\to r$ induced by the counit $rl\to \id_y$. 
\end{proposition}
\begin{proof}
The representable $2$-functor $\mathfrak{D}(-,y): \mathfrak{D}^{1\op}\to \cat$ induces an adjunction of categories 
	\begin{equation}
		\begin{tikzcd}
			\mathfrak{D}(y,y):\ar[r,shift left=2,"r^*"{name=A}]&\mathfrak{D}(x,y)\ar[l,shift left=2,"l^*"{name=B}]\ar[phantom,from=A,to=B,"\dashv"rotate=-90]
		\end{tikzcd}
	\end{equation}
	In particular, for all $f\in \mathfrak{D}(x,y)$, $g\in \mathfrak{D}(y,y)$, the counit map $rl\to \id_y$ induces an equivalence 
	\begin{equation}
		\mathfrak{D}(x,y)(gr,f)\simeq \mathfrak{D}(x,y)(r^*g,f)\simeq \mathfrak{D}(y,y)(g,l^*f)\simeq \mathfrak{D}(y,y)(g,fl).
	\end{equation}
	Putting $f=r$, for all $g\in \mathfrak{D}(x,y)$, the counit map induces an equivalence $\mathfrak{D}(x,y)(gr,r)\simeq \mathfrak{D}(y,y)(g,rl)$, exhibiting $rl$ as an endomorphism object for $r$. 	
\end{proof}

Fix a monoidal category $\mathcal{C}$.  For each $\mathcal{C}$-module category $\mathcal{M}$ there is a forgetful functor $\Lmod(\mathcal{M})\to \alg(\mathcal{C})$, taking a module object in $\mathcal{M}$ to the algebra which acts on it.  The fiber over an algebra $A$ is the category $\Lmod_A( \mathcal{M})$.  We would like to show that the assignment 
\begin{equation}
	\mathcal{M}\mapsto (\Lmod(\mathcal{M})\to \alg(\mathcal{C}))
\end{equation}
is (lax monoidally) functorial in $\mathcal{M}\in \Lmod_ \mathcal{C}$.  For this, we shall introduce one hypothesis

\begin{hypothesis}\label{hyp:MonoidalCoCart}
	Let $\mathcal{O}$ be an operad and $\mathcal{C}$ an $\mathcal{O}$-monoidal $2$-category. There is an equivalence of categories of $2$-functors
	\begin{equation}
		\alg_ \mathcal{O}\fun(\mathcal{C},\cat)\simeq \mathrm{cocart}^\mathcal{O}(\mathcal{C}). 
	\end{equation}
	In particular, using the universal property of Day convolution, every $\mathcal{O}$-monoidal coCartesian fibration determines a lax $\mathcal{O}$-monoidal $2$-functor $\mathcal{C}\to \cat$. 
\end{hypothesis}
The $1$-categorical version of this is proven in \cite{Hinich2015} Proposition A.2.1, \cite{RamziUnstraightening} Theorem 2.1.  This version is stronger in that it gives an equivalence of categories of $\cat$-enriched functors, rather than just categories of $\mathcal{S}$-enriched functors.

\begin{construction}\label{LModunfold}
	Let $\mathcal{K}$ be a collection of categories.  Suppose that $\mathcal{C}\in \alg_{E_1\otimes \mathcal{O}}(\cat(\mathcal{K}))$.  By Proposition \ref{prop:AlgMonoidal}, there is $\mathcal{O}$-monoidal category $\Lmod_\mathcal{C}(\cat(\mathcal{K}))^\otimes\to \mathcal{O}^\otimes$.  The universal property of the Boardman-Vogt tensor product induces a functor $\LM^\otimes \otimes \Lmod_\mathcal{C}(\cat(\mathcal{K}))^\otimes\to \cat(\mathcal{K})^\otimes$.  This unstraightens to a coCartesian fibration, which we denote by 
\begin{equation}
	\widetilde{\Lmod_\mathcal{C}(\cat(\mathcal{K}))}^\otimes \to LM^\otimes \otimes \Lmod_\mathcal{C}(\cat(\mathcal{K}))^\otimes. 
\end{equation}
Passing to $\LM$- and $E_1$-algebras and using Proposition \ref{prop:AlgMonoidal}, we obtain coCartesian fibrations 
\begin{align}
p:\alg_{\LM}(\widetilde{\Lmod_\mathcal{C}(\cat(\mathcal{K}))})\to \Lmod_\mathcal{C}(\cat(\mathcal{K}))^\otimes\\
q:\alg(\widetilde{\Lmod_\mathcal{C}(\cat(\mathcal{K})})\to \Lmod_\mathcal{C}(\cat(\mathcal{K}))^\otimes
\end{align}
Pulling back along the map $E_1\to LM$ induces a morphism $\phi$ of coCartesian fibrations over $\Lmod_ \mathcal{C}(\cat(\mathcal{K}))^\otimes$. 
\begin{equation}
	\begin{tikzcd}[column sep=5]
		\alg_{LM}(\widetilde{\Lmod_\mathcal{C}(\cat(\mathcal{K}))})^\otimes\ar[rr,"\phi"]\ar[dr,swap,"p"]&&\alg(\widetilde{\Lmod_\mathcal{C}(\cat(\mathcal{K}))})^\otimes\ar[dl,"q"]\\
											     &\Lmod_\mathcal{C}(\cat(\mathcal{K}))^\otimes
	\end{tikzcd}
\end{equation}
Note that each fiber is an object in $\widehat{\cat}(\mathcal{K})$ by \cite{LurieHA} 4.2.3.5.  Then by Hypothesis \ref{hyp:MonoidalCoCart}, straightening over $\Lmod_\mathcal{C}(\cat(\mathcal{K}))$ gives a pair of lax $\mathcal{O}$-monoidal 2-functors and a monoidal natural transformation between them 
\begin{equation}
	\begin{tikzcd}[column sep=80]
		\Lmod_\mathcal{C}(\cat(\mathcal{K})) 
		\ar[bend left=30,"\alg_{LM}",name=U]{r}[name=U]{}
		\ar[bend right=30,swap,"\alg"]{r}[name=D]{} &
		       \widehat{\cat}(\mathcal{K})
		       \arrow[shorten <=10pt,shorten >=10pt,Rightarrow,to path={(U) -- node[label=right:$\mathrm{st}(\phi)$]{} (D)}]
			    \end{tikzcd}
	\end{equation}
This gives a 2-functor $\Theta:\Lmod_\mathcal{C}(\cat(\mathcal{K}))\to \fun([1],\widehat{\cat}(\mathcal{K}))$. 
\end{construction}

\noindent Suppose that $T\in \alg \Omega\mathfrak{D}$ and $x\in \mathfrak{D}$.  Our goal is to  show that the assignment
\begin{equation}
	(T,x)\mapsto \Lmod_T(\mathfrak{D}(x,\mathds{1}))
\end{equation}
defines a lax monoidal functor which is contravariant in $T$. 
\begin{construction}\label{EMRepresentable}
	Suppose that $\mathcal{K}$ is a collection of categories and $\mathcal{C}\in \alg_\mathcal{O}(\cat(\mathcal{K}))$.  Consider the lax $\mathcal{O}$-monoidal $2$-functor $\Theta:\Lmod_\mathcal{C}(\cat(\mathcal{K}))\to \fun([1],\widehat{\cat}(\mathcal{K}))$ defined in Construction \ref{LModunfold}, which takes $\mathcal{M}\in \Lmod_\mathcal{C}(\cat(\mathcal{K}))$ to $\mathrm{st}(\phi)(\mathcal{M}):\alg_{\LM}(\mathcal{M})\to \alg(\mathcal{C})$.  For all $\mathcal{M}\in \Lmod_\mathcal{C}(\cat(\mathcal{K}))$, $\mathrm{st}(\phi)(\mathcal{M})$ is a Cartesian fibration, and morphisms in $\Lmod_\mathcal{C}(\cat(\mathcal{K}))$ are taken to morphisms of Cartesian fibrations by \cite[Corollary 4.2.3.2]{LurieHA}. 
	Then $\Theta$ factors $\mathcal{O}$-monoidally through the category  $\mathrm{cart}(\alg(\mathcal{C}))$ of Cartesian fibrations over $\alg(\mathcal{C})$.  By straightening over $\alg(\mathcal{C})$ using Hypothesis \ref{hyp:MonoidalCoCart}, we obtain lax $\mathcal{O}$-monoidal $2$-functor 
	\begin{equation}
		\Lmod^*_{-}(-):	\Lmod_\mathcal{C}(\cat(\mathcal{K}))\to \fun(\alg(\mathcal{C})^{\op},\cat(\mathcal{K}))^\otimes. 
	\end{equation}
One should note that, for algebra morphism $f:A\to B$ in $\alg(\mathcal{C})$ and $\mathcal{M}\in \Lmod_\mathcal{C}(\cat(\mathcal{K}))$, the functor $\Lmod_B(\mathcal{M})\to \Lmod_A(\mathcal{M})$ is given by restriction of scalars along $f$. 
Specializing Construction \ref{Yonedamodule} to $\mathcal{V}=\cat$, and composing with $\Lmod^*_-(-)$, we obtain lax $\mathcal{O}$-monoidal $2$-functor 
	\begin{equation}
		\mathfrak{D}^{1\op}\to \Lmod_{\Omega\mathfrak{D}}(\cat)\to \fun(\alg(\Omega\mathfrak{D})^{\op},\cat(\mathcal{K}))^\otimes,
	\end{equation}
which takes $x\in \mathfrak{D}$ to the functor $\Lmod_{-}(\mathfrak{D}(x,\mathds{1})):\alg(\Omega\mathfrak{D})^{\op}\to \cat(\mathcal{K})$.  Using $\mathcal{O}$-monoidal tensor hom adjunction, we have $\mathcal{O}$-monoidal equivalence
\begin{multline}
	\fun(\mathfrak{D}^{1\op},\fun(\Omega\alg(\mathfrak{D})^{\op},\cat(\mathcal{K})))^\otimes\simeq \fun(\mathfrak{D}^{1\op}\otimes \alg(\Omega\mathfrak{D})^{\op},\cat(\mathcal{K}))^\otimes\\
	\simeq \fun(\alg(\Omega\mathfrak{D})^{\op},\fun(\mathfrak{D}^{1\op},\cat(\mathcal{K})))^\otimes,
\end{multline}
which takes the previously defined lax $\mathcal{O}$-monoidal functor to lax $\mathcal{O}$-monoidal functor 
\begin{equation}
	\Lmod^*_{-}(\mathfrak{D}(-,\mathds{1})):	\alg(\Omega\mathfrak{D})^{\op}\to \fun(\mathfrak{D}^{1\op},\cat(\mathcal{K}))^\otimes,
\end{equation}
with algebra morphism $f:A\to B$ in $\alg(\Omega\mathfrak{D})$ taken to restriction of scalars along $f$.
\end{construction}

\begin{proposition}[\cite{StefanichPhD} 7.3.7]\label{Stefanich7.3.7}
	Suppose that $\mathfrak{D}$ is a $2$-category, $y\in \mathfrak{D}$, and $T\in \alg(\mathfrak{D}(y,y))$ is a monad on $y$.  Suppose that $g:X\to y$ is a module for $T$, and $g$ exhibits $X$ as an Eilenberg-Moore object for $T$.  Then $g$ admits a right adjoint and the induced morphism $T\to \End(g)$ is an equivalence. 
\end{proposition}

\begin{notation}
	Suppose that $\mathfrak{D}$ is an $2$-category and $y\in \mathfrak{D}$.  Denote by $\mathfrak{D}^R_{/y}$ the full subcategory of $\mathfrak{D}_{/y}$ on the right adjoint morphisms to $y$. 
\end{notation}

\begin{theorem}\label{EMR}
	Suppose that $\mathfrak{D}\in \alg_\mathcal{O}(2\cat)$ admits Eilenberg-Moore objects at $\mathds{1}\in \mathfrak{D}$.  Then, there is an adjunction 
	\begin{equation}
		\begin{tikzcd}
			\alg(\Omega \mathcal{D})^{\op}\ar[r,shift left=2,"\EM_R"{name=B}]&		\mathfrak{D}_{/\mathds{1}}^R \ar[l,shift left=2,"\End_L"{name=A}] \ar[phantom,from=A,to=B,"\dashv"rotate=90] 
		\end{tikzcd}
	\end{equation}
$\EM_R$ is fully faithful and lax $\mathcal{O}$-monoidal, taking a monad $T$ on $\mathds{1}$ to its Eilenberg-Moore object, and morphisms of monads to restriction of scalars, and $\End_L$ takes a morphism to its endomorphism monad. 
\end{theorem}

\begin{proof}
	$\Lmod^*_{-}(\mathfrak{D}(-,\mathds{1})):\alg(\Omega \mathcal{D})^{\op}\to \fun(\mathfrak{D}^{1\op},\cat)$ factors monoidally through the full subcategory on representable presheaves, by the assumption that $\mathfrak{D}$ admits Eilenberg-Moore objects.  
	Since $\alg(\Omega \mathfrak{D})$ has initial object given by the algebra structure on $\id_{\mathds{1}}$, which is taken to $\mathfrak{D}(-,\mathds{1})$ under the above functor, the desired functor is given by 
\begin{equation}
	\EM_R:\alg(\Omega\mathfrak{D})^{\op}\simeq (\alg(\Omega\mathfrak{D})_{\id/})^{\op}\to \fun(\mathfrak{D}^{1\op},\cat)^{rep}_{/\mathfrak{D}(-,\mathds{1})}\simeq \mathfrak{D}_{/\mathds{1}},
\end{equation}
where the final equivalence is by monoidal Yoneda embedding, and factors through $\mathfrak{D}_{/\mathds{1}}^R$ since $\mathfrak{D}$ admits Eilenberg-Moore objects at $\mathds{1}$.   

Every object in $\mathfrak{D}_{/\mathds{1}}^R$ admits an endomorphism object in $\Omega\mathfrak{D}$ by Proposition \ref{adjendmon}.  To show that $\EM_R$ admits a left adjoint, it suffices to show that there is a natural equivalence 
\begin{equation}
\alg(\Omega\mathcal{D})(T, \mathrm{End}(g))\simeq \alg(\Omega\mathcal{D})^\op(\mathrm{End}(g),T)\simeq \mathcal{D}^R_{/\mathds{1}}(g, \EM_R(T))
\end{equation}
for all $g\in \mathfrak{D}^R_{/\mathds{1}}$.  This is immediate from the definition of Eilenberg-Moore object and endomorphism monad, since both functors universally encode composition actions of $T$ on $g$.  

Finally, to prove that $\mathrm{EM}_R$ is fully faithful, it suffices to show that the counit of this adjunction is an equivalence.  But this is exactly Proposition \ref{Stefanich7.3.7}.

\end{proof}
\begin{remark}
	A global (non-monoidal) version of Theorem \ref{EMR} for structured monads appears in \cite{HeineMonadicity} Theorem 6.30.  Functoriality of $\End_L:\mathfrak{D}_{/y}^R\to \alg(\mathfrak{D}(y,y))^\op$ treated in \cite{HeineMonadicity} 6.28. 

\end{remark}

\begin{corollary}
	Suppose that $\mathfrak{D}\in \alg_\mathcal{O}(2\cat)$ and $T\in \alg_{E_1\otimes \mathcal{O}}(\Omega\mathfrak{D})$.  Then $\EM(T)\in \mathrm{co}\alg_\mathcal{O}(\mathfrak{D})$. 
\end{corollary}

\subsection{The left Eilenberg-Moore functor}

Suppose that $T\in \alg \Omega\mathfrak{D}$.  Our goal is to  show that the assignment
\begin{equation}
	(T,x)\mapsto \Lmod_T(\mathfrak{D}(x,\mathds{1}))
\end{equation}
defines a lax monoidal functor which is \textbf{covariant} in $T$.

\begin{construction}\label{EMop}
	Let $\mathcal{K}$ be a collection of categories containing $\Delta^\op$.  Suppose that $\mathcal{C}\in \alg_{E_1\otimes\mathcal{O}}\cat(\mathcal{K})$.  Consider the $2$-functor $\Theta:\Lmod_\mathcal{C}(\cat(\mathcal{K}))\to \fun([1],\widehat{\cat}(\mathcal{K}))$ given in Construction \ref{LModunfold}, which takes $\mathcal{M}\in \Lmod_\mathcal{C}(\cat(\mathcal{K}))$ to $\mathrm{st}(\phi)(\mathcal{M}):\alg_{\LM}(\mathcal{M})\to \alg(\mathcal{C})$. 

	Since $\mathcal{K}$-contains $\Delta^{\op}$, $\mathrm{st}(\phi)(\mathcal{M})$ is a coCartesian fibration by \cite[Proposition 4.6.2.17]{LurieHA}.  Moreover, an edge $(A,L)\to (B,N)$ in $\alg_{\LM}(\mathcal{M})$ is coCartesian over $\alg(\mathcal{C})$ if and only if it induces an equivalence $B\otimes_AL\simeq N$.  Since $\Theta$ takes a morphism in $\Lmod_\mathcal{C}(\cat(\mathcal{K}))$ to a natural transformation which preserves geometric realizations, we see that morphisms in $\Lmod_\mathcal{C}(\cat(\mathcal{K}))$ are taken to morphisms of coCartesian fibrations.  Then $\Theta$ factors $\mathcal{O}$-monoidally through $\mathrm{cocart}(\alg(\mathcal{C}))$, the category of coCartesian fibrations over $\alg(\mathcal{C})$.  By straightening over $\alg(\mathcal{C})$ and using Hypothesis \ref{hyp:MonoidalCoCart}, we obtain lax $\mathcal{O}$-monoidal $2$-functor 
	\begin{equation}
		\Lmod_{!-}(-):\Lmod_\mathcal{C}(\cat(\mathcal{K}))\to \fun(\alg(\mathcal{C}),\cat(\mathcal{K}))^\otimes. 
	\end{equation}
One should note that, for algebra morphism $f:A\to B$ in $\alg(\mathcal{C})$ and $\mathcal{M}\in \Lmod_\mathcal{C}(\cat(\mathcal{K}))$, the functor $\Lmod_A(\mathcal{M})\to \Lmod_B(\mathcal{M})$ is given by extension of scalars along $f$. 	
	Composing with the lax $\mathcal{O}$-monoidal functor of Constructions \ref{Yonedamodule}, applied to $\mathcal{V}=\cat(\mathcal{K})$, gives a lax $\mathcal{O}$-monoidal $2$-functor 
	\begin{equation}
		\mathfrak{D}^{1\op}\to \Lmod_{\Omega\mathfrak{D}}(\cat(\mathcal{K}))\to \fun(\alg(\Omega\mathfrak{D}),\cat(\mathcal{K}))^\otimes. 
	\end{equation}
Using $\mathcal{O}$-monoidal tensor hom adjunction we have $\mathcal{O}$-monoidal equivalence
\begin{multline}
	\fun(\mathfrak{D}^{1\op},\fun(\alg(\Omega\mathfrak{D}),\cat(\mathcal{K})))^\otimes\simeq \fun(\mathfrak{D}^{1\op}\otimes \alg(\Omega\mathfrak{D}),\cat(\mathcal{K}))^\otimes\\
	\simeq \fun(\alg(\Omega\mathfrak{D}),\fun(\mathfrak{D}^{1\op},\cat(\mathcal{K})))^\otimes,
\end{multline}
which takes the previously defined lax $\mathcal{O}$-monoidal $2$-functor to a lax $\mathcal{O}$-monoidal $2$-functor
\begin{equation}
	\Lmod_{!-}(\mathfrak{D}(-,\mathds{1})):\alg(\Omega\mathfrak{D})\to \fun(\mathfrak{D}^{1\op},\cat(\mathcal{K}))^\otimes, 
\end{equation}
with algebra morphism $f:A\to B$ in $\alg(\Omega\mathfrak{D})$ taking to extension of scalars along $f$. 
\end{construction}

\begin{notation}
 Suppose that $\mathfrak{D}$ is an $2$-category and $y\in \mathfrak{D}$.  Denote by $\mathfrak{D}_{y/}^L$ the full subcategory of $\mathfrak{D}_{/y}$ on the left adjoint morphisms from $y$. 
\end{notation}

\begin{theorem}\label{EML}
	Let $\mathcal{K}$ be a collection of categories containing $\Delta^{\op}$, and suppose that $\mathfrak{D}\in \alg_\mathcal{O}(2\cat(\coprod^{\geq 1} \mathcal{K})$ is a $2$-category which is locally compatible with $\mathcal{K}$ indexed colimits, and also admits Eilenberg-Moore objects at $\mathds{1}\in \mathfrak{D}$.  Then there is lax $\mathcal{O}$-monoidal functor 
	\begin{equation}
\EM_L:\alg(\Omega\mathfrak{D})\to \mathfrak{D}_{\mathds{1}/}^L
	\end{equation}
	taking a monad $T$ to the left-adjoint monadic morphism $\mathds{1}\to \EM(T)$, and morphisms of monads to extension of scalars.  
\end{theorem}
\begin{proof}
	$\Lmod_{!-}(\mathfrak{D}(-,\mathds{1})):\alg(\Omega\mathfrak{D})\to \fun(\mathfrak{D}^{1\op},\cat(\mathcal{K}))$ factors $\mathcal{O}$-monoidally through the full subcategory on representable presheaves by the assumption that $\mathfrak{D}$ admits Eilenberg-Moore objects.  $\alg(\Omega\mathfrak{D})$ has an initial object $id_ \mathds{1}$.  The desired functor is given by 
	\begin{equation}
\EM_L:		\alg(\Omega\mathfrak{D})\simeq \alg(\Omega\mathfrak{D})_{id_ \mathds{1}/}\to \fun(\mathfrak{D}^{1\op},\cat)_{\mathfrak{D}(\_,\mathds{1})/}^{rep}\simeq \mathfrak{D}_{\mathds{1}/}.
	\end{equation}
\end{proof}

\begin{remark}
	When $\mathfrak{D}$ does not admit all Eilenberg-Moore objects at $\mathds{1}$, we still get a functor of the form described above by restricting to the full subcategory of $\alg(\Omega \mathfrak{D})$ on those monads which do admit Eilenberg-Moore objects. 
\end{remark}

\begin{remark}
	\cite{HaugsengMonad} corollary 5.8 provides a closely related result for `cosmifiable' $2$-categories.  \cite{ReutterZetto} section 4 provides a similar result for $\pr^\mathcal{V}$. 
\end{remark}

\begin{corollary}
	Suppose that $\mathfrak{D}\in \alg_\mathcal{O}(2\cat)$ locally admits $\Delta^\op$-indexed colimits and $T\in \alg_{E_1\otimes \mathcal{O}}(\Omega \mathfrak{D})$.  Then $\EM(T)\in \alg_\mathcal{O}(\mathfrak{D})$. 
\end{corollary}

\subsection{Application to relative tensor products of modules}

\begin{construction}
	Suppose that $\mathcal{V}\in \alg_{E_2\otimes \mathcal{O}}(\pr)$, $\mathcal{K}$ is a collection of small categories, and $\mathcal{C}\in \alg_{E_1\otimes \mathcal{O}}(\cat^\mathcal{V}(\mathcal{K}))$.  Then $\Rmod_\mathcal{C}(\cat^\mathcal{V}(\mathcal{K}))\in \alg_\mathcal{O}(\pr^{\cat^\mathcal{V}(\mathcal{K})})$ by \cite{LurieHA} 4.8.5.16 applied to a larger universe.  We may then consider $\Rmod_\mathcal{C}(\cat^\mathcal{V}(\mathcal{K}))$ as a large $2$-category which is locally compatible with $\mathcal{K}$-indexed colimits by changing enrichment along the underlying category functor $\cat^\mathcal{V}(\mathcal{K})\to \cat(\mathcal{K})$.
\end{construction}
\begin{construction}\label{algmonad}
	Suppose that $\mathcal{V}\in \alg_{E_1\otimes \mathcal{O}}(\pr)$, $\mathcal{K}$ is a collection of categories containing $\Delta^\op$, and $\mathcal{C}\in \alg_{\mathcal{O}}(\cat^\mathcal{V}(\mathcal{K}))$.  There is $\mathcal{O}$-monoidal equivalence $\Rmod_\mathcal{C}(\cat^\mathcal{V}(\mathcal{K}))\simeq \fun^\mathcal{V}(B\mathcal{C}^\op,\cat^\mathcal{V}(\mathcal{K}))$.  By delooping hypothesis (Theorem \ref{delooping}), we get $E_1\otimes \mathcal{O}$-monoidal equivalence $\mathcal{C}\simeq \Omega_ \mathcal{C}\Rmod_\mathcal{C}(\cat^\mathcal{V}(\mathcal{K}))$.  In particular, we get an $\mathcal{O}$-monoidal equivalence 
	\begin{align}
		\alg(\mathcal{C})\simeq \alg(\Omega_ \mathcal{C}\Rmod_\mathcal{C}\cat^\mathcal{V}(\mathcal{K}))&&A\mapsto A\otimes -
	\end{align}
	between algebra objects in $\mathcal{C}$ and right $\mathcal{C}$-linear monads on $\mathcal{C}$.  Note that we have an equivalence of categories $\Lmod_{A}(\mathcal{C})\simeq \Lmod_{A\otimes-}(\mathcal{C})$. 
\end{construction}

\begin{corollary}\label{AlgebraEMR}
	Suppose that $\mathcal{V}\in \alg_{E_2\otimes \mathcal{O}}(\pr)$, $\mathcal{K}$ is a small collection of small categories, and $\mathcal{C}\in \alg_{E_1\otimes \mathcal{O}}(\cat^\mathcal{V}(\mathcal{K}))$.  Then there is a fully faithful $\mathcal{O}$-monoidal functor 
	\begin{equation}
		\Lmod^*_{-}(\mathcal{C}):\alg(\mathcal{C})^\op\to \Rmod_\mathcal{C}(\cat^\mathcal{V}(\mathcal{K}))_{/\mathcal{C}}
	\end{equation}
such that for $A\in \alg(\mathcal{C})$, the induced morphism $\Lmod_A(\mathcal{C})\to \mathcal{C}$ is the right adjoint in a monadic adjunction in $\Rmod_\mathcal{C}(\cat^\mathcal{V}(\mathcal{K}))$.  Moreover, $\Lmod^*_{-}(\mathcal{C})$ admits a left adjoint. 
\end{corollary}
\begin{proof}
$\Rmod_\mathcal{C}(\cat^\mathcal{V})$ admits monoidal Eilenberg-Moore objects since it admits $\cat$-weighted limits.  Then the desired functor is given by 
\begin{equation}
	\begin{tikzcd}
		\alg(\mathcal{C})^\op\ar[r,"\ref{algmonad}"]&\Omega_ \mathcal{C}\Rmod_\mathcal{C}(\cat^\mathcal{V})^\op\ar[r,"\EM_R"]&\Rmod_\mathcal{C}(\cat^\mathcal{V})_{/\mathcal{C}}.
	\end{tikzcd}
\end{equation}

\end{proof}

\begin{corollary}\label{AlgebraEML}
	Suppose that $\mathcal{V}\in \alg_{E_2\otimes \mathcal{O}}(\pr)$, $\mathcal{K}$ is a collection of categories containing $\Delta^\op$, and $\mathcal{C}\in \alg_{E_1\otimes \mathcal{O}}(\cat^\mathcal{V}(\mathcal{K}))$.  Then there is a $\mathcal{O}$-monoidal functor 
	\begin{equation}
		\Lmod_{!-}(\mathcal{C}):\alg(\mathcal{C})\to \Rmod_\mathcal{C}(\cat^\mathcal{V}(\mathcal{K}))_{\mathcal{C}/}
	\end{equation}
such that for $A\in \alg(\mathcal{C})$, the induced morphism $\mathcal{C}\to \Lmod_A(\mathcal{C})$ is the left adjoint in a monadic adjunction in $\Rmod_\mathcal{C}(\cat^\mathcal{V}(\mathcal{K}))$.  
\end{corollary}
\begin{proof}
$\Rmod_\mathcal{C}(\cat^\mathcal{V}(\mathcal{K}))$ is locally compatible with geometric realizations, and admits Eilenberg-Moore objects since it admits $\cat$-weighted limits.  Then the desired functor is given by 
\begin{equation}
	\begin{tikzcd}
		\alg(\mathcal{C})\ar[r,"\ref{algmonad}"]&\Omega_ \mathcal{C}\Rmod_\mathcal{C}(\cat^\mathcal{V})\ar[r,"\EM_L"]&\Rmod_\mathcal{C}(\cat^\mathcal{V}(\mathcal{K}))_{\mathcal{C}/}.
	\end{tikzcd}
\end{equation}
\end{proof}

\begin{example}
	Corollaries \ref{AlgebraEMR} and \ref{AlgebraEML} appear to be very similar.  One should notice that \ref{AlgebraEML} requires that $\mathcal{C}$ is is locally compatible with geometric realizations, while \ref{AlgebraEMR} does not.  It is useful to provide an example to illustrate the difference.
	\begin{itemize}
		\item If one has bialgebra object $A$ in $\mathcal{C}$, (i.e. $A\in co\alg(\alg(\mathcal{C}))$, the monoidality of $\Lmod_{-}^*(\mathcal{C})$ gives the natural monoidal structure on $\Lmod_A(\mathcal{C})$. 
		\item Alternatively, if one gives $E_{k+1}$ algebra object $A$ in $\mathcal{C}$, monoidality of $\Lmod_{!-}(\mathcal{C})$ gives the natural $E_k$-monoidal structure on $\Lmod_A(\mathcal{C})$. 
	\end{itemize}

\end{example}

We use an alternate proof to obtain that the analogous functor for presentably $\mathcal{V}$-enriched categories is also fully faithful with right adjoint.  We are not sure to what extent this applies beyond the presentable case. 

\begin{lemma}\label{presentableAlgEML}
	Suppose that $\mathcal{V}\in \alg_{E_{k+1}}(\pr)$ for $k\geq 1$, and $\mathcal{C}\in \alg_{E_{k}}(\pr^\mathcal{V})$.  Then there is a fully faithful $E_{k-1}$-monoidal functor 
	\begin{equation}
		\Rmod_{!-}(\mathcal{C}):\alg(\mathcal{C})\to \Lmod_\mathcal{C}(\pr^\mathcal{V})_{\mathcal{C}/}
	\end{equation}
such that for $A\in \alg(\mathcal{C})$, the induced morphism $\mathcal{C}\to \Lmod_A(\mathcal{C})$ is the left adjoint in a monadic adjunction in $\Rmod_\mathcal{C}(\pr^\mathcal{V})$.  

Moreover, $\Lmod_{!-}(\mathcal{C})$ admits a right adjoint. 
\end{lemma}

\begin{proof}
	It follows from corollary \ref{AlgebraEML} (or alternatively \cite{LurieHA} 4.8.5.16) applied to a larger universe that one obtains an $E_{i-1}$-monoidal functor $\Lmod_{!-}(\mathcal{C}):\alg(\mathcal{C})\to \Rmod_\mathcal{C}(\pr)_{\mathcal{C}/}$.  By \cite{LurieHA}, this is fully faithful with right adjoint.  By \cite{LurieHA} Theorem 7.1.3.1 and corollary 7.1.3.4, one has a monoidal equivalence $\Rmod_\mathcal{C}(\pr^\mathcal{V})\simeq \Rmod_{\mathcal{C}}(\pr)$.  The desired functor is given by the composition
	\begin{equation}
		\alg(\mathcal{C})\xrightarrow{\Lmod_{!-}(\mathcal{C})}\Rmod_{\mathcal{C}}(\pr)_{/\mathcal{C}}\simeq \Rmod_\mathcal{C}(\pr^\mathcal{V}).
	\end{equation}
	
\end{proof}

\begin{remark}
	If $\mathcal{C}$ is a monoidal category which is compatible with $\Delta^\op$-indexed colimits, then one could give a definition of the Morita category as in \cite{LurieHA} 4.8.4.9 and \cite{SoergelBimodule} section 6 as the essential image of the functor $\Lmod_{!-}(\mathcal{C}):\alg(\mathcal{C})\to \Rmod_ \mathcal{C}(\cat^\mathcal{V})$.  This inherits a natural enrichment over $\cat^\mathcal{V}(\mathcal{K})$ from the enrichment of $\Rmod_\mathcal{C}(\cat^\mathcal{V}(\mathcal{K}))$.  It follows from \cite{LurieHA} Theorem 4.8.4.1 that the morphisms in this category agree with the (enriched) category of $B-A$ bimodule objects in $\mathcal{C}$. 

	This provides an enriched version of the `higher Morita category' and `even higher Morita catgory' constructions of \cite{ScheimbauerPhD,Haugseng2017,JohnsonFreydScheimbauer,Haugseng2017,Haugseng2023}.  Note however, that this Morita category does not generally admit geometric realizations, so we may not be able to take iterated Morita categories in this way.
\end{remark}

\begin{warning}
	The version of Morita category given above uses the opposite convention for left/right modules to others listed above.  This is fixed by instead considering the essential image of $\Rmod_{!-}(\mathcal{C}):\alg(\mathcal{C})\to \Lmod_\mathcal{C}(\pr^\mathcal{V})$. 
\end{warning}

\printbibliography
\vspace{0.5cm}
\begin{flushleft}
Centre for Quantum Mathematics,
University of Southern Denmark\\
\textit{Email address}:  stockall@imada.sdu.dk\\
\textit{Website}: devonstockall.com
\end{flushleft}
\end{document}